\documentclass[a4paper,12pt]{article}
\makeatletter
    
    \@addtoreset{equation}{subsection}
  \makeatother
\usepackage{amsmath,amssymb,amsthm}
\usepackage{amscd}
\usepackage[all]{xy}
\theoremstyle{definition}
\newtheorem{dfn}{Definition}[subsection]
\newtheorem{exple}[dfn]{Example}

\theoremstyle{plain}
\newtheorem{prop}[dfn]{Proposition}
\newtheorem{thm}[dfn]{Theorem}
\newtheorem{lem}[dfn]{Lemma}
\newtheorem{cor}[dfn]{Corollary}

\begin{document}
\title{The Meyer functions for projective varieties and their application to local signatures for
fibered 4-manifolds}
\author{Yusuke Kuno}
\date{}
\maketitle

\begin{abstract}
We study a secondary invariant, called the Meyer function, on the fundamental
group of the complement of the dual variety of a smooth projective variety.
This invariant have played an important role when studying the local signatures of
fibered 4-manifolds from topological point of view. As an application of our
study, we define a local signature for generic non-hyperelliptic fibrations
of genus 4 and 5 and compute some examples. 
\end{abstract}

\section{Introduction}
Let $\Sigma_g$ be a closed oriented $C^{\infty}$-surface of genus $g\ge 0$.
\textit{The mapping class group of} $\Sigma_g$, which we denote by $\Gamma_g$,
is the group of orientation preserving diffeomorphisms of $\Sigma_g$ modulo isotopy.
The group cohomology of $\Gamma_g$ attracts attentions because: 1)
its element plays a characteristic class of oriented $\Sigma_g$-bundles, 2)
over the rational coefficients, it is isomorphic to the cohomology of the moduli
space of compact Riemann surfaces of genus $g$.

As for the degree two part, the cohomology group itself has been determined. Harer \cite{Harer} proved that
$H^2(\Gamma_g;\mathbb{Z})\cong \mathbb{Z}$ for $g\ge 3$; $H^2(\Gamma_1;\mathbb{Z})\cong \mathbb{Z}/12\mathbb{Z}$
and $H^2(\Gamma_2;\mathbb{Z})\cong \mathbb{Z}/10\mathbb{Z}$ are classically known. However,
as a reflection of the fact that $\Gamma_g$ is related to various mathematical objects,
there have been known various 2-cocycles of $\Gamma_g$ arising from different contexts.

One of these is \textit{Meyer's signature cocycle $\tau_g$}, introduced by
W.\ Meyer \cite{Meyer} and rediscovered later by Turaev \cite{Turaev}.
The definition involves the signature of 4-manifolds and will be recalled
in this section.

The main object we study here is the \textit{Meyer function}, a secondary invariant 
associated to $\tau_g$. The work of Meyer \cite{Meyer} is considered as the origin of it.
He showed that: for $g=1$ or $2$, there exists a unique $\mathbb{Q}$-valued 1-cochain
$\phi_g\colon \Gamma_g \rightarrow \mathbb{Q}$ whose coboundary equals to $\tau_g$.
He also gave an explicit formula for $\phi_1$. Note that $\Gamma_1$ is isomorphic to
$SL(2;\mathbb{Z})$. Atiyah \cite{Atiyah} showed interesting aspects of $\phi_1$
as a function: he showed that the value of $\phi_1$
for a hyperbolic element $\alpha\in SL(2;\mathbb{Z})$
coincides with various values associated to $\alpha$ such as the special value of
a Shimizu $L$-function determined by $\alpha$, an arithmetic invariant,
or the $\eta$-invariant of the mapping torus of $\alpha$, a differential geometric invariant.

Recently there are several works that give higher genera or higher dimensional analogues
of $\phi_1$ or $\phi_2$. In \cite{Endo,Mor} the Meyer function on the hyperelliptic mapping
class group is studied. In \cite{Endo}, application to the local signature for hyperelliptic
fibrations is dealt and in \cite{Mor}, a relation to the $\eta$-invariant of mapping tori
is studied. In \cite{Iida}, as a higher dimensional generalization of $\phi_2$,
the Meyer function for the family of smooth theta divisors is studied.

\noindent \textbf{The Meyer function for a projective variety.}
In this paper we give other analogues of Meyer's $\phi_1$ or $\phi_2$ and
discuss their applications.
We consider the family of Riemann surfaces constructed as follows.

Let $X\subset \mathbb{P}_N$ be a smooth complex projective variety of dimension $n$, embedded
in the complex projective space of dimension $N$. Throughout the paper, we assume that $$N>n\ge 2.$$
Let $k:=N-n+1$. We denote by $G_k(\mathbb{P}_N)$ the set of all $k$-planes of $\mathbb{P}_N$.
Let $D_X$ be the set of $k$-planes meeting $X$ not transversally.
When $n=2$, $D_X$ is the classical dual variety of $X$.
In \cite{GKZ}, $D_X$ is called the $k$-th associated subvariety of $X$.
Let $U^X=G_k(\mathbb{P}_N)\setminus D_X$.
For $W\in U^X$, $W$ and $X$ meet transversally so their intersection $X\cap W$ has a natural structure
of a compact Riemann surface. Thus setting
$$\mathcal{C}^X:=\left\{ (x,W)\in \mathbb{P}_N \times U^X; x\in X\cap W \right\},$$
the second projection
\begin{equation}
\label{eq:1-1}
p_X:=p_2|_{\mathcal{C}^X} \colon \mathcal{C}^X \rightarrow U^X
\end{equation}
is a complex analytic family of compact Riemann surfaces. Let $g$ be the genus of the fibers and
$$\rho_X\colon \pi_1(U^X)\rightarrow \Gamma_g$$
the topological monodromy (see the conventions below) of (\ref{eq:1-1}).
Let $\rho_X^*\tau_g$ be the pull back of $\tau_g$ by $\rho_X$.

\begin{thm}[$=$ Theorem \ref{thm:3-1-1}]
\label{thm:1-0-1}
There exists a unique $\mathbb{Q}$-valued 1-cochain $\phi_X\colon \pi_1(U^X)\rightarrow \mathbb{Q}$
whose coboundary equals to $\rho_X^*\tau_g$.
In particular, we have $\rho_X^*[\tau_g]=0\in H^2(\pi_1(U^X);\mathbb{Q})$.
\end{thm}

We remark here that if $\rho\colon G \rightarrow \Gamma_g$ is a homomorphism from a group $G$
to $\Gamma_g$ and $\phi\colon G\rightarrow \mathbb{Q}$ is an 1-cochain whose coboundary equals to
$\rho^*\tau_g$, then $\phi$ is always a class function: $\phi(xyx^{-1})=\phi(y)$ for
$x,y\in G$. This is easily derived from properties of $\tau_g$ (see \cite{Ku}, Appendix).
In particular the above $\phi_X$ is a class function.
We call $\phi_X$ \textit{the Meyer function associated to $X\subset \mathbb{P}_N$}.
In fact this theorem is regarded as a further generalization of \cite{Ku},
where the case of $X$ being the $d$-th Veronese image of $\mathbb{P}_2$ is studied.
Our proof is based on a geometric feature of $\tau_g$, and applications of
the Novikov additivity of the signature are essential.

There are several studies on the fundamental group of the complement of the dual variety
(or more generally the associated subvariety in a Grassmannian), for example, see \cite{CT,DL}
and a recent work of I.\ Shimada \cite{Shi}. However, it is still a mysterious object and lots
of the properties are unknown.
The function $\phi_X$ tells us some information as to $\pi_1(U^X)$
under a mild condition, see Proposition \ref{prop:3-6-1}.

\noindent \textbf{An application to local signatures.}
One reason to seek a generalization of Meyer's $\phi_1$ or $\phi_2$ comes from the motivation
to treat \textit{localization of the signature} for the case of the genus of the fibers
greater than two via Meyer functions.

Let us go back to the case of $g=1$ or $2$ for explanation.
The coboundary condition $\delta \phi_g=\tau_g$ leads to
an immediate consequence: for an oriented surface bundle of genus $\le 2$
over a closed oriented surface, the signature of the total space is zero.
Proceeding further, let $M$ (resp. $B$) be a closed oriented
$C^{\infty}$-manifold of dimension 4 (resp. 2) and $f\colon M\rightarrow B$ a proper surjective
$C^{\infty}$-map having a structure of surface bundle of genus $g$, over the outside of finitely
many points $b_1,\ldots,b_m \in B$. We call such a triple $(M,f,B)$ a \textit{fibered 4-manifold}.
The fiber germ $\mathcal{F}_i$ over $b_i$ is called a singular fiber germ.
Typical examples are elliptic surfaces or Lefschetz fibrations.

In the above situation, the advantage of $\phi_1$ or $\phi_2$ is that we can associate each singular fiber
germ with its local invariant $\sigma(\mathcal{F}_i)\in \mathbb{Q}$, called the \textit{local signature}.
The adjective "local" comes from the equality
$${\rm Sign}(M)=\sum_{i=1}^m\sigma(\mathcal{F}_i).$$
The definition of $\sigma(\mathcal{F}_i)$ is given by
\begin{equation}
\label{eq:1-2}
\sigma(\mathcal{F}_i)=\phi_g(x_i)+{\rm Sign}(N(f^{-1}(b_i))),
\end{equation}
where $x_i \in \Gamma_g$ is the local monodromy around $b_i$ and ${\rm Sign}(N(f^{-1}(b_i)))$ is the
signature of a fiber neighborhood of $f^{-1}(b_i)$. This formulation first appeared in 
Y. Matsumoto's papers \cite{Mat1,Mat2}.
For generalizations of this story for higher genera, there is an obstruction:
the class $[\tau_g]$ is a generator of $H^2(\Gamma_g;\mathbb{Q})\cong \mathbb{Q}$ for $g\ge 3$.

Local signatures are also studied from complex geometric or algebro geometric point of view, see \cite{AK,AY}.
In these setting, a local signature is defined by another way and can be defined even if
$g\ge 3$, by assigning some algebro geometric conditions on the general fibers.
There is an important point to note here: when $g\ge 3$,
there is a fiber germ with a non-trivial local signature
but topologically being a trivial $\Sigma_g$-bundle. To capture such phenomena,
it is insufficient just looking at the shape of $f^{-1}(b_i)$ or the local monodromy $x_i$,
hence we need to modify (\ref{eq:1-2}).

As for higher genera analogues of Y. Matsumoto's approach,
Endo \cite{Endo} studied the local signature for hyperelliptic fibrations.
In \cite{Ku} non-hyperelliptic fibrations of genus 3 are discussed.
In this paper using the Meyer functions $\phi_X$ for particular choices of $X$,
we will discuss non-hyperelliptic fibrations of genus 4 or 5.
The modification of (\ref{eq:1-2}) is achieved by introducing a group with some
universal property and the Meyer function on that group. $x_i$ in (\ref{eq:1-2}) is
replaced by the \textit{lifted monodromy}, see Definition \ref{dfn:4-1-8}.
One advantage of our local signature is that we only need the complex structures on the general
fibers so it is not necessary $f$ itself should be holomorphic.
Although we don't know whether our local signature is the same as the others \cite{AK,AY},
we will observe the coincidence on some examples of singular fiber germs.

\noindent \textbf{Organization of the paper.} Section 2 is a preparation for section 3.
We describe the tangent space of $D_X$ and study the situation when a holomorphic disk
intersects $D_X$ transversally. These considerations will be used in Proposition \ref{prop:3-2-1}.
In section 3 we prove Theorem \ref{thm:1-0-1} by a purely topological argument.
Using the method of Lefschetz pencils, we give a formula for the value
of $\phi_X$ on a special element, called \textit{lasso}. We also study the second bounded cohomology of $\pi_1(U^X)$.
In section 4 applications to non-hyperelliptic fibrations of genus 4 or 5 are discussed.
First we explain our approach to local signatures via Meyer functions (see Proposition \ref{prop:4-1-7}),
then proceed to the particular cases. When the genus is 4, we assign the general fibers to be "of rank 4".
When the genus is 5, we assign the general fibers to be non-trigonal.
In the case of genus 4, we compute the value of our local signature
for some fiber germs.

\vspace{0.5cm}
In the rest of this introduction we fix conventions and recall Meyer's signature
cocycle.

\noindent \textbf{Topological monodromy.} We adopt the following:
1) for any two mapping classes $f_1$ and $f_2$, the multiplication $f_1\circ f_2$ means
that $f_2$ is applied first,
2) for any two homotopy classes of based loops $\ell_1$ and $\ell_2$, their product
$\ell_1\cdot \ell_2$ means that $\ell_1$ is traversed first.

Let $p\colon \mathcal{C}\rightarrow B$ be an oriented $\Sigma_g$-bundle. Choose a base point $b_0\in B$ and
fix an identification $\phi\colon \Sigma_g \stackrel{\cong}{\rightarrow} p^{-1}(b_0)$.
For each based loop $\ell\colon [0,1]\rightarrow B$
the pull back $\ell^*\mathcal{C}\rightarrow [0,1]$ by $\ell$ is a trivial $\Sigma_g$-bundle.
Hence there exists a trivialization
$\Phi\colon \Sigma_g \times [0,1]\rightarrow \ell^*\mathcal{C}$ such that $\Phi(x,0)=\phi(x)$, $x\in \Sigma_g$.
By assigning the isotopy class of $\Phi(\cdot,1)^{-1}\circ \phi$ to the homotopy class of $\ell$,
we obtain a map $\rho$, called the \textit{topological monodromy of
$p\colon \mathcal{C}\rightarrow B$}, from $\pi_1(B,b_0)$ to $\Gamma_g$.
Under the conventions above, $\rho$ is a homomorphism.

\noindent \textbf{Meyer's signature cocycle.}
Let $P$ denote the pair of pants, i.e., $P=S^2\setminus \bigcup_{i=1}^{3}{\rm Int}D_i$
where $D_i$, $1\le i\le 3$, are the three disjoint closed disks in the 2-sphere $S^2$.
Choose a base point $p_0\in {\rm Int}P$ and fix based loops $\ell_i$, $1\le i\le 3$
such that each $\ell_i$ is homotopic to the loop traveling once the boundary $\partial D_i$
by counter clockwise manner, and the product $\ell_1\cdot \ell_2\cdot \ell_3$
is null homotopic.
For $(f_1,f_2)\in \Gamma_g\times \Gamma_g$, we can construct an oriented $\Sigma_g$-bundle
$E(f_1,f_2)$ over $P$ such that the topological monodromy
$\pi_1(P,p_0)\rightarrow \Gamma_g$ sends $[\ell_i]$ to $f_i$ for $i=1,2$.
$E(f_1,f_2)$ is a compact $C^{\infty}$-manifold of dimension 4 endowed
with the natural orientation. Thus the signature of $E(f_1,f_2)$
is defined and we set $$\tau_g(f_1,f_2):={\rm Sign}(E(f_1,f_2)).$$
By the Novikov additivity of the signature $\tau_g$ turns out to be a 2-cocycle of $\Gamma_g$.
The class $[\tau_g]\in H^2(\Gamma_g;\mathbb{Z})$ equals to
$1/3$ times the first MMM class \cite{Miller, Morita, Mumford}.

There is a linear algebraic description of $\tau_g$ given in \cite{Meyer}.
Let $\Gamma_g\rightarrow Sp(2g;\mathbb{Z})$ be the homomorphism obtained by 
the action of $\Gamma_g$ on the first homology of $\Sigma_g$, and let $A_1$ and $A_2$
be the image of $f_1$ and $f_2$ by this homomorphism, respectively.
Let
$$J=\left( \begin{array}{cc} 0 & I_g \\ -I_g & 0 \\ \end{array} \right)$$
where $I_g$ is the $g\times g$ identity matrix, and consider the linear space 
$$V_{A_1,A_2}:=\{ (x,y)\in \mathbb{R}^{2g}\oplus \mathbb{R}^{2g}; (A_1^{-1}-I_{2g})x+(A_2-I_{2g})y=0\},$$
where $I_{2g}$ is the $2g\times 2g$ identity matrix. Then
$$\langle (x,y),(x^{\prime},y^{\prime}) \rangle_{A_1,A_2}=\ ^t(x+y)J(I_{2g}-A_2)y^{\prime}$$
turns out to be a symmetric bilinear form on $V_{A_1,A_2}$ hence its signature is defined.
As proved in \cite{Meyer}, we have
\begin{equation}
\label{eq:1-3}
\tau_g(f_1,f_2)={\rm Sign}(V_{A_1,A_2},\langle \ ,\ \rangle_{A_1,A_2}).
\end{equation}

Here we correct some errors about signs in \cite{Ku}. In Appendix of \cite{Ku},
we have adopted the same notations about topological monodromies as this paper and
have defined $\tau_g(f_1,f_2)=-{\rm Sign}(E(f_1,f_2))$.
Then Definition 7.1(p.\ 943) should be corrected as
$${\rm loc.sig}^{\mathcal{Q}}(\mathcal{F}):=-\phi^4(\theta(\mathcal{F}^0)_*(\gamma))+{\rm Sign}(E).$$
The equation ${\rm Sign}(\pi^{-1}(B_0))=\sum_{i=1}^n\phi^4(\theta(\mathcal{F}_i^0)_*(\gamma))$
in the proof of Theorem 7.2(p.\ 944) should be corrected by multiplying the right hand side by $-1$.
The proof of Proposition 5.1(p.\ 936) should be corrected similarly and all the values of the Meyer function
in \cite{Ku} should be multiplied by $-1$.

\noindent \textbf{Notations.}
For an 1-cochain $\phi\colon G\rightarrow A$ of a group $G$ with coefficient in an abelian group $A$,
its coboundary is meant the map $\delta \phi\colon G\times G\rightarrow A$ defined by
$$\delta \phi(x,y)=\phi(x)-\phi(xy)+\phi(y).$$

For a complex manifold $M$, we denote by $K_M$ the canonical divisor of $M$. 
More generally, for a possibly singular variety $Y$, we denote by $\omega_Y$ the dualizing sheaf of $Y$.
We use this notion only when $Y$ is given as a hypersurface in a complex manifold $M$. In this case
$\omega_Y$ is an invertible sheaf on $Y$ given by the adjunction formula:
$$\omega_Y=(K_M+Y)|_Y.$$ 

For integers $p,q$ with $0<p<q$, we denote by $G_{p,q}$
the Grassmannian of all $p$-planes of $\mathbb{C}^q$. Note that $G_k(\mathbb{P}_N)$
is naturally isomorphic to $G_{k+1,N+1}$.

\section{Preliminaries from complex algebraic geometry}
In this section we describe some properties of $D_X$.
When $n=2$, $D_X$ is an irreducible variety in
$G_{N-1}(\mathbb{P}_N)=\mathbb{P}_N^{\vee}$, the dual projective space of $\mathbb{P}_N$, and $D_X$ is
classically known as \textit{the dual variety of $X$}. In fact, the treatment here is a generalization of the
treatments in sections 1 and 2 of K.\ Lamotke's paper \cite{La} to the case of general $n$. Corollary
\ref{cor:2-2-3}, Proposition \ref{prop:2-3-3}, and Theorem \ref{thm:2-3-4} will be used in later sections.
Let
$$\mathcal{W}:=\left\{ (x,W)\in \mathbb{P}_N \times G_k(\mathbb{P}_N); x\in X\cap W \right\}.$$
Then there are two projections
$p_1\colon \mathcal{W}\rightarrow X$,
and
$$p_2\colon \mathcal{W}\rightarrow G_k(\mathbb{P}_N).$$

\subsection{Coordinate description of $p_2$}
In the following we give an explicit coordinate description of $p_2$.

Let $(x^0,W_0)\in \mathcal{W}$.
By choosing appropriate homogeneous coordinates $[x_0:x_1:\cdots :x_N]$ of $\mathbb{P}_N$,
we may assume that $x^0=[1:0:\cdots:0]$ and $W_0$ is given by $x_{k+1}=\cdots =x_N=0$.

We first introduce local coordinates of $\mathcal{W}$ near $(x^0,W_0)$.
For $x\in X$, $p_1^{-1}(x)$ is the set
of $k$-planes of $\mathbb{P}_N$ through $x$, which is isomorphic to $G_{k,N}$.
The open set $\{ x_0\neq 0\}$ of $\mathbb{P}_N$ is identified with $\mathbb{C}^N$ by
\begin{equation}
\label{eq:2-1-1}
[x_0:x_1:\cdots :x_N]\mapsto \left(x_1/x_0,\ldots,x_N/x_0 \right).
\end{equation}
Thus for $(x,\stackrel{\circ}{W})\in (X\cap \{ x_0\neq 0\})\times G_{k,N}$, considering the affine subspace
$x+\stackrel{\circ}{W}\subset \mathbb{C}^N\cong \{ x_0\neq 0\}$ and taking its closure in $\mathbb{P}_N$,
we have a trivialization of $p_1$ over $X\cap \{ x_0\neq 0\}$:
\begin{equation}
\label{eq:2-1-2}
(X\cap \{ x_0\neq 0\})\times G_{k,N}\stackrel{\cong}{\longrightarrow}
p_1^{-1}(X\cap \{ x_0\neq 0\}).
\end{equation}
Let $e_i:=(0,\ldots,\stackrel{i}{1},\ldots,0)^t\in \mathbb{C}^N$ for $1\le i\le N$. In view of
(\ref{eq:2-1-1}), $\stackrel{\circ}{W}_0:=W_0 \cap \{ x_0\neq 0\}$ is the $k$-plane of
$\mathbb{C}^N$ spanned by $e_1,\ldots,e_k$. Let $W_1$ be the subspace of $\mathbb{C}^N$
spanned by $e_{k+1},\ldots,e_N$.
For $F\in {\rm Hom}(\stackrel{\circ}{W}_0,W_1)$ let
$\stackrel{\circ}{W}(F):={\rm Span}(f_1,\ldots,f_k)\in G_{k,N}$,
where $f_i:=e_i+F(e_i)$, $1\le i\le k$. By this mapping $F\mapsto \stackrel{\circ}{W}(F)$,
$\mathcal{U}:=\{ \stackrel{\circ}{W}\in G_{k,N}; \stackrel{\circ}{W}\cap W_1=0 \}$ is identified with
${\rm Hom}(\stackrel{\circ}{W}_0,W_1)$.
Introducing $\{ F_i^j\}$ by
$$F(e_i)=\sum_{j=k+1}^NF_i^je_j,\ 1\le i\le k$$
for $F\in {\rm Hom}(\stackrel{\circ}{W}_0,W_1)$,
the set of functions $\{ F_i^j \}$ serves as local coordinates of $\mathcal{U}$.

Choose a sufficiently small local coordinate neighborhood $(U;t_1,\ldots,t_n)$ of $X$ centered at $x^0$
(i.e., $x^0$ corresponds to the origin $(0,\ldots,0)$) so that $U\subset X\cap \{ x_0\neq 0\}\subset \mathbb{C}^N$.
Then points in $U$ can be expressed as $x(t_1,\ldots,t_n)=(x_1,\ldots,x_N)$ where
$x_i=x_i(t_1,\ldots,t_n)$, $1\le i\le N$ are holomorphic in $t_1,\ldots,t_n$.

In view of (\ref{eq:2-1-2}),
\begin{equation}
\label{eq:2-1-3}
\left( U\times \mathcal{U};t_1,\ldots,t_n,\{ F_i^j \} \right)
\end{equation}
can be used as a local coordinate neighborhood of $\mathcal{W}$ centered at $(x^0,W_0)$.

Next we introduce local coordinates of $G_k(\mathbb{P}_N)$ near $W_0$.
Let $\hat{e}_i=(0,\ldots,\stackrel{i+1}{1},\ldots,0)^t \in \mathbb{C}^{N+1}$, 
$0\le i\le N$ and $\hat{W}_0={\rm Span}(\hat{e}_0,\hat{e}_1,\ldots,\hat{e}_k)$,
$\hat{W}_1={\rm Span}(\hat{e}_{k+1},\ldots,\hat{e}_{N})$.
Then by the natural isomorphism
$G_k(\mathbb{P}_N)\cong G_{k+1,N+1}$, $W_0$ corresponds to $\hat{W}_0$ (recall that $W_0$ is given by
$x_{k+1}=\cdots x_N=0$).
For $G\in {\rm Hom}(\hat{W}_0,\hat{W}_1)$ let $W(G):={\rm Span}(g_0,g_1,\ldots,g_k)$,
where $g_i:=\hat{e}_i+G(\hat{e}_i)$, $0\le i\le k$.
By this mapping $G\mapsto W(G)$, $\mathcal{V}:=\{ W\in G_{k+1,N+1}; W\cap \hat{W}_1=0\}$ is
identified with ${\rm Hom}(\hat{W}_0,\hat{W}_1)$.
Introducing $\{G_i^j \}$ by
$$G(\hat{e}_i)=\sum_{j=k+1}^NG_i^j \hat{e}_j,\ 0\le i\le k.$$
for $G\in {\rm Hom}(\hat{W}_0,\hat{W}_1)$, then
\begin{equation}
\label{eq:2-1-4}
\left(\mathcal{V},\{ G_i^j \} \right)
\end{equation}
is a local coordinate neighborhood of $G_{k+1,N+1}\cong G_k(\mathbb{P}_N)$ near $W_0$.

Now for $(t_1,\ldots,t_n,\{ F_i^j \})\in U\times \mathcal{U}$, using the local coordinates
(\ref{eq:2-1-3}) and (\ref{eq:2-1-4}), write
$p_2(t_1,\ldots,t_n,\{ F_i^j \})=G(t_1,\ldots,t_n,\{F_i^j \})\in {\rm Hom}(\hat{W}_0,\hat{W}_1)$.
The closure of the affine space $x(t_1,\ldots,t_n)+\stackrel{\circ}{W}(F)\subset \mathbb{C}^N$
in $\mathbb{P}_N$ must be equal to $W(G(t_1,\ldots,t_n,\{ F_i^j \}))$.
On the other hand, at the level of $G_{k+1,N+1}$,
the closure is a $(k+1)$-plane spanned by $(1,x_1,\ldots,x_n)^t$ and $\hat{f}_1,\ldots,\hat{f}_k$,
where for $1\le i\le k$, $\hat{f}_i\in \mathbb{C}^{N+1}$ is the image of $f_i$ under the inclusion
$\mathbb{C}^N \hookrightarrow \mathbb{C}^{N+1}$, $(y_1,\ldots,y_N)\mapsto (0,y_1,\ldots,y_N)$.
From these we establish the following:
\begin{lem}[local description of $p_2$]
\label{lem:2-1-1}
Let $( U\times \mathcal{U};t_1,\ldots,t_n,\{ F_i^j \} )$ and
$(\mathcal{V},\{ G_i^j \} )$ be the local coordinate neighborhoods of $\mathcal{W}$ and
$G_k(\mathbb{P}_N)$ respectively, as above. Let
$G_i^j=G_i^j(t_1,\ldots,t_n,\{ F_i^j \})$ be the local coordinates of the point
$W(G(t_1,\ldots,t_n,\{ F_i^j \}))$. Then we have
$$G_i^j=F_i^j \textrm{ for } 1\le i\le k,\ k+1\le j\le N, \textrm{ and}$$
\begin{equation}
\label{eq:2-1-5}
G_0^j=x_j-\sum_{i=1}^kx_iF_i^j, \textrm{ for } k+1\le j\le N.
\end{equation}
\end{lem}

\subsection{Irreducibility of $D_X$}
The $k$-th associated subvariety $D_X$ is irreducible. Although this might be
well known, here we include the proof of it together with the irreducibility
of some loci related to $D_X$.

\begin{dfn}
\label{dfn:2-2-1}
Define the subsets of $\mathcal{W}$ as follows:
$$\mathcal{D}:=\left\{ (x,W)\in \mathcal{W}; T_xX+T_xW\neq T_x\mathbb{P}_N \right\},$$
and for integers $i\ge 1$,
\begin{eqnarray*}
\mathcal{Y}_i &:=& \left\{ (x,W)\in \mathcal{W}; \dim(T_xX+T_xW)=N+1-i \right\} \\
&=& \left\{ (x,W)\in \mathcal{W}; \dim(T_xX \cap T_xW)=i \right\}.
\end{eqnarray*}
Finally, for $i\ge 1$, $\mathcal{D}_i:=\bigcup_{j\ge i}\mathcal{Y}_j$.
\end{dfn}

Note that $\mathcal{Y}_i$ is empty for $i> \max(n,k)$ and $\mathcal{Y}_i=\mathcal{D}_i \setminus \mathcal{D}_{i+1}$.
Also we have $p_2(\mathcal{D})=D_X$. All $\mathcal{D}_i$ are closed analytic subsets of $\mathcal{W}$ and
we have a filtration
$$\mathcal{D}_1=\mathcal{W} \supset \mathcal{D}_2=\mathcal{D} \supset \mathcal{D}_3 \supset
\mathcal{D}_4 \supset \cdots \supset \mathcal{D}_{\max(n,k)}.$$

Using Lemma \ref{lem:2-1-1}, we can verify that
$\mathcal{D}$ is the set of critical points of $p_2\colon \mathcal{W}\rightarrow G_k(\mathbb{P}_N)$
and $\mathcal{Y}_i$ (resp.\ $\mathcal{D}_i$) is the set of points $(x,W)\in \mathcal{W}$ such that
the differential $(p_2)_* \colon T_{(x,W)}\mathcal{W}\rightarrow T_WG_k(\mathbb{P}_N)$
has corank $i-1$ (resp.\ corank$\ge i-1$).

\begin{thm}
\label{thm:2-2-2}
\begin{enumerate}
\item
For each $i$, $\mathcal{Y}_i$ is a connected submanifold of $\mathcal{W}$ with codimension $i^2-i$
and is open and dense in $\mathcal{D}_i$.
\item
For each $i$, $\mathcal{D}_i$ is an irreducible analytic subset of $\mathcal{W}$ with codimension $i^2-i$.
\end{enumerate}
\end{thm}

Since $\dim \mathcal{W}=\dim G_k(\mathbb{P}_N)+1$ and $p_2$ is a proper
holomorphic map, we get the following which we will use later.

\begin{cor}
\label{cor:2-2-3}
The set $D_X$ is an irreducible analytic subset of $G_k(\mathbb{P}_N)$ with codimension $\ge 1$.
If the codimension of $D_X$ is 1, $p_2(\mathcal{D}_3)$ is a proper analytic subset of $D_X$.
\end{cor}

\begin{proof}[Proof of Theorem \ref{thm:2-2-2}]
Let us introduce some notations. Let $V$ be a fixed $n$-dimensional
subspace of $\mathbb{C}^N$ and for $i=1,2,\ldots$, let
$$Y_i:=\left\{ W\in G_{k,N}; \dim(V+W)=N+1-i \right\},$$
and $D_i=\bigcup_{j\ge i}Y_j$.
By trivializations of the pair of holomorphic vector bundles $(T\mathbb{P}_N|_X,TX)$ on $X$,
we see that $p_1|_{\mathcal{D}_i}$ and $p_1|_{\mathcal{Y}_i}$ are holomorphic fiber bundle with fiber
isomorphic to $D_i$ and $Y_i$, respectively. Since $Y_i$ is open and dense in $D_i$, $\mathcal{Y}_i$
is also open and dense in $\mathcal{D}_i$.
We see that $Y_i$ has a structure of a connected complex manifold of dimension
$kn-k+i-i^2$. This can be seen by considering the projection
$Y_i\rightarrow G_i(V)$, $W\mapsto W\cap V$ ($G_i(V)$ is the Grassmannian
of all $i$-planes of $V$). This shows $\mathcal{Y}_i$ is a connected
complex manifold with the desired codimension. 

We next prove the second part. Since $\mathcal{Y}_i$ is contained in $\mathcal{D}_i$, the first part shows
that the set of smooth points of $\mathcal{D}_i$ is connected, hence $\mathcal{D}_i$ is irreducible.
Also we have $\dim \mathcal{D}_i=\dim \mathcal{Y}_i$. This completes the proof.
\end{proof}

\subsection{The tangent space of $D_X$ for the case $D_X$ is a hypersurface}
In this subsection we describe the tangent space of $D_X$ at a generic point of $D_X$ under the assumption
that the codimension of $D_X$ is 1, i.e., $D_X$ is a hypersurface of $G_k(\mathbb{P}_N)$. Then
$p_2|_{\mathcal{D}}\colon \mathcal{D}\rightarrow D_X$ is a dominant regular map between projective varieties
of the same dimension. By Sard's lemma for varieties (see Chapter 3 of \cite{Mum} for instance), there exists
a proper analytic subset $E^{\prime}\subset D_X$ such that
\begin{enumerate}
\item $E^{\prime}$ contains $S(D_X)$, the set of singular points of $D_X$,
\item the differential $(p_2|_{\mathcal{D}})_*\colon T_{(x,W)}\mathcal{D}\rightarrow T_WD_X$ is an
isomorphism for $(x,W)\in (p_2|_{\mathcal{D}})^{-1}(D_X\setminus E^{\prime})\setminus S(\mathcal{D})$,
where $S(\mathcal{D})$ denotes the set of singular points of $\mathcal{D}$.
\end{enumerate}
By Theorem \ref{thm:2-2-2}, $\mathcal{Y}_2$ is contained in $\mathcal{D}\setminus S(\mathcal{D})$
so $S(\mathcal{D})\subset \mathcal{D}_3$. Setting $E:=E^{\prime}\cup p_2(\mathcal{D}_3)$
then this is a proper analytic subset of $D_X$ by Corollary \ref{cor:2-2-3}. Now we have

\begin{lem}
\label{lem:2-3-1}
Suppose that the codimension of $D_X$ is 1. Then there exists a proper analytic subset $E$ of $D_X$ such that
\begin{enumerate}
\item $E$ contains $S(D_X)$,
\item $(p_2|_{\mathcal{D}})^{-1}(D_X\setminus E)\subset \mathcal{Y}_2$, in particular
 $(p_2|_{\mathcal{D}})^{-1}(D_X\setminus E)$ is contained in $\mathcal{D}\setminus S(\mathcal{D})$,
\item the differential $(p_2|_{\mathcal{D}})_*$is an
isomorphism at $(x,W)\in (p_2|_{\mathcal{D}})^{-1}(D_X\setminus E)$.
\end{enumerate}
\end{lem}

Suppose $(x^0,W_0)\in \mathcal{D}$ and $W_0=p_2(x^0,W_0)\in D_X\setminus E$. By Lemma \ref{lem:2-3-1},
$(x^0,W_0)\in \mathcal{Y}_2$. Then the following relations among the subspaces of $T_{W_0}G_k(\mathbb{P}_N)$ hold:
$$(p_2)_*(T_{(x^0,W_0)}\mathcal{W})\supset (p_2)_*(T_{(x^0,W_0)}\mathcal{Y}_2)
=(p_2|_{\mathcal{D}})_*(T_{(x^0,W_0)}\mathcal{D})=T_{W_0}D_X.$$
Since $(x^0,W_0)\in \mathcal{Y}_2$, $(p_2)_*(T_{(x^0,W_0)}\mathcal{W})$ has codimension 1 in
$T_{W_0}G_k(\mathbb{P}_N)$ (see subsection 2.2). Therefore, we have
\begin{equation}
\label{eq:2-3-1}
T_{W_0}D_X=(p_2)_*(T_{(x^0,W_0)}\mathcal{W}).
\end{equation}
Recall the local coordinates of $\mathcal{W}$ and $G_k(\mathbb{P}_N)$ in subsection 2.1. 
By Lemma \ref{lem:2-1-1}, we see that $(p_2)_*(T_{(x^0,W_0)}\mathcal{W})$ is generated by the $n$ vectors
\begin{equation}
\label{eq:2-3-2}
\frac{\partial x_{k+1}}{\partial t_i}(0)\frac{\partial}{\partial G_0^{k+1}}+
\ldots +\frac{\partial x_{N}}{\partial t_i}(0)\frac{\partial}{\partial G_0^N},
\end{equation}
$1\le i\le n$, and the $k(n-1)$ vectors
\begin{equation}
\label{eq:2-3-3}
\frac{\partial}{\partial G_i^j},
\end{equation}
$1\le i\le k$, $k+1\le j\le N$. Let
$$\frac{\partial \boldsymbol{x}}{\partial t_i}=\left(\frac{\partial x_{k+1}}{\partial 
t_i},\ldots,\frac{\partial x_{N}}{\partial t_i}\right)^t\in \mathbb{C}^{n-1},\ 1\le i\le n.$$
Since $(x^0,W_0)\in \mathcal{Y}_2$ the rank of the matrix
$$\left(\displaystyle\frac{\partial \boldsymbol{x}}{\partial t_1}(0),\ldots,
\displaystyle\frac{\partial \boldsymbol{x}}{\partial t_n}(0) \right)$$
is $(n-2)$. By an arrangement of indices we may assume that
\begin{equation}
\label{eq:2-3-4}
\frac{\partial \boldsymbol{x}}{\partial t_i}(0)
,\ 1\le i\le n-2,\ \textrm{ are linearly independent,}
\end{equation}
then the other two column vectors are in the linear span of these.

\begin{prop}
\label{prop:2-3-2}
Suppose that the codimension of $D_X$ is 1 and let $E$ be as in Lemma \ref{lem:2-3-1}. Let $(x^0,W_0)\in \mathcal{D}$
and suppose $W_0\in D_X\setminus E$. Then under the assumption (\ref{eq:2-3-4}),
the tangent space $T_{W_0}D_X$ is given by
\begin{equation}
\label{eq:2-3-5}
T_{W_0}D_X=\left\{ \sum_{i,j}u_i^j \frac{\partial}{\partial G_i^j};
\det \left( \boldsymbol{u},
\frac{\partial \boldsymbol{x}}{\partial t_1}(0),
\ldots,
\frac{\partial \boldsymbol{x}}{\partial t_{n-2}}(0) \right)=0 \right\}.
\end{equation}
Here, $\boldsymbol{u}=(u_0^{k+1},\ldots,u_0^{N})^t\in \mathbb{C}^{n-1}$.
\end{prop}
\begin{proof}
By the assumption (\ref{eq:2-3-4}), the right hand side of (\ref{eq:2-3-5}) is a hyperplane of 
$T_{W_0}G_k(\mathbb{P}_N)$. Also, the vectors (\ref{eq:2-3-2}) and (\ref{eq:2-3-3}) are clearly contained
in the right hand side of (\ref{eq:2-3-5}). This completes the proof.
\end{proof}

Regarding $p_2\colon \mathcal{W}\rightarrow G_k(\mathbb{P_N})$ as a family of algebraic curves,
we investigate its pull back by a mapping into $G_k(\mathbb{P_N})$ which does not meet $E$ and is
transverse to $D_X$. We first show that the total space of the pull back has the structure of a manifold.

\begin{prop}
\label{prop:2-3-3}
Suppose that the codimension of $D_X$ is 1 and let $E$ be as in Lemma \ref{lem:2-3-1}.
Let $B$ be a $C^{\infty}$-manifold of dimension $\ge 2$ and let $\iota\colon B\rightarrow G_k(\mathbb{P}_N)$
be a $C^{\infty}$-map satisfying $\iota^{-1}(E)=\emptyset$ and transverse to $D_X$. Then, the pull back
\begin{eqnarray*}
\iota^*\mathcal{W}&:=& \left\{ (b,(x,W))\in B\times\mathcal{W}; \iota(b)=W \right\} \\
&\cong & \left\{ (b,x)\in B \times X; x\in \iota(b) \right\}
\end{eqnarray*}
of $p_2$ by $\iota$ has the natural structure of a $C^{\infty}$-manifold as a $C^{\infty}$-submanifold of
$B\times X$. Moreover, if $B$ is a complex manifold and $\iota$ is a holomorphic map, $\iota^*\mathcal{W}$
has the natural structure of a complex manifold as a complex submanifold of $B\times X$.
\end{prop}
\begin{proof}
We only treat the case $B$ is the small disk $\Delta:=\{ z\in \mathbb{C};|z|<\varepsilon\},\ \varepsilon >0$
and $\iota$ is a holomorphic map such that $\iota^{-1}(D_X)=\{ 0\}$.
A similar argument proves the general case (see also Lemma 2.4 in \cite{Ku}).

By the assumption, we have $\iota(0)\in D_X\setminus E$ and the transversality
\begin{equation}
\label{eq:2-3-6}
\iota_*(T_0\Delta)+T_{\iota(0)}D_X=T_{\iota(0)}G_k(\mathbb{P}_N).
\end{equation}
Let $(z_0,x^0)\in \iota^*\mathcal{W}$ and write $W_0:=\iota(z_0)$. Choosing the local coordinates of
$\mathcal{W}$ and $G_k(\mathbb{P}_N)$ as in subsection 2.1 we denote by $\iota_i^j$ the coordinate
expression of $\iota$ with respect to the local coordinates $\{ G_i^j\}$. In particular,
we have $\iota_i^j(z_0)=0$. By Lemma \ref{lem:2-1-1},
the local equation of $\iota^*\mathcal{W}$ near $(z_0,x^0)$ is given by
$$-\iota_0^j(z)+x_{j}(t_1,\ldots,t_n)-\sum_{i=1}^kx_i(t_1,\ldots,t_n)\iota_i^j(z)=0$$
for $k+1\le j \le N$. Let $\psi_j(z,t_1,\ldots,t_n)$ be the left hand side of the above equation.
The Jacobian matrix of $(\psi_{k+1},\ldots,\psi_N)$ at $(z_0,0,\ldots,0)$ is the $(n-1)\times (n+1)$ matrix
\begin{equation}
\label{eq:2-3-7}
\left( -\boldsymbol{\iota^{\prime}}(z_0),
\frac{\partial \boldsymbol{x}}{\partial t_1}(0),\ldots,
\frac{\partial \boldsymbol{x}}{\partial t_n}(0) \right),
\end{equation}
where $\boldsymbol{\iota^{\prime}}=(d{\iota_0^{k+1}}/dz,\ldots,d{\iota_0^{N}}/dz)$.
We claim that this matrix is of full rank. Suppose $(x^0,W_0)\notin \mathcal{D}$, namely $(x^0,W_0)\in \mathcal{Y}_1$.
Then the $(n-1)\times n$ matrix obtained by deleting the first column of (\ref{eq:2-3-7})
is already of full rank so is (\ref{eq:2-3-7}).

Suppose $(x^0,W_0)\in \mathcal{D}$. Then by the assumption, $z_0=0$ and $W_0\in D_X\setminus E$.
Proposition \ref{prop:2-3-2} and (\ref{eq:2-3-6}) shows that
\begin{equation}
\label{eq:2-3-8}
\det \left( \boldsymbol{\iota^{\prime}}(0),
\frac{\partial \boldsymbol{x}}{\partial t_1}(0),\ldots,
\frac{\partial \boldsymbol{x}}{\partial t_{n-2}}(0) \right)\neq 0,
\end{equation}
therefore (\ref{eq:2-3-7}) is of full rank also in this case.
By the implicit function theorem the assertion follows.
\end{proof}

Let $\Delta$ and $\iota$ be as in the proof of Proposition \ref{prop:2-3-3}. The pull back
$\iota^*\mathcal{W}=\{ (z,x)\in \Delta \times X; x\in \iota(z) \}$ has the natural projection
$f_{\iota}\colon \iota^*\mathcal{W}\rightarrow \Delta$. Explicitly, $f_{\iota}$ is given by
$f_{\iota}(z,x)=z$.

\begin{thm}
\label{thm:2-3-4}
Notations are as above. Then, $(z_0,x^0)\in \iota^*\mathcal{W}$
is a critical point of $f_{\iota}$ if and only if $z_0=0$
and $(x^0,\iota(0))\in \mathcal{D}$. All the critical points are non-degenerate.
\end{thm}
In fact, we will see in Corollary \ref{cor:3-4-4} that there is only one critical point.
By an argument like the Morse lemma, we see that near each critical point $f_{\iota}$ looks like
$(z_1,z_2)\mapsto z_1^2+z_2^2$. The next subsection is devoted to the proof of this theorem.

\subsection{Proof of Theorem \ref{thm:2-3-4}}
We will use the notations in the proof of Proposition \ref{prop:2-3-3}.

Let $(z_0,x^0)\in \iota^*\mathcal{W}$ and write $W_0:=\iota(z_0)$. Suppose $(x^0,W_0)\notin \mathcal{D}$.
Then we may assume that in the column vectors of (\ref{eq:2-3-7})
$$\frac{\partial \boldsymbol{x}}{\partial t_1}(0),\ldots,
\frac{\partial \boldsymbol{x}}{\partial t_{n-1}}(0)$$
are linearly independent. By the implicit function theorem, there exist local coordinates $(s_1,s_2)$ of
$\iota^*\mathcal{W}$ centered at $(z_0,x^0)$ such that the points near $(z_0,x^0)$ can be expressed
as $(z,t_1,\ldots,t_{n-1},t_n)$ where
$$z=s_1,t_1=t_1(s_1,s_2),\ldots,t_{n-1}=t_{n-1}(s_1,s_2),t_n=s_2,$$
and $t_i(s_1,s_2)$ are holomorphic in $s_1,s_2$. Since $f_{\iota}(z,x)=z=s_1$,
$(z_0,x^0)$ is not a critical point of $f_{\iota}$.

Suppose $(x^0,W_0)\in \mathcal{D}$. Then as we have seen in the proof of Proposition \ref{prop:2-3-3}, $z_0=0$,
$W_0\in D_X \setminus E$, and $(x^0,W_0)\in \mathcal{Y}_2$. We may assume (\ref{eq:2-3-4}). Then we have
the inequality (\ref{eq:2-3-8}) hence there exist local coordinates $(s_1,s_2)$ of $\iota^*\mathcal{W}$
centered at $(z_0,x^0)$ such that
\begin{equation}
\label{eq:2-4-1}
z=z(s_1,s_2),t_1=t_1(s_1,s_2),\ldots,t_{n-2}=t_{n-2}(s_1,s_2),t_{n-1}=s_1,t_{n}=s_2.
\end{equation}
For $k+1\le j\le N$, differentiating the identity
$$\psi_j(z(s_1,s_2),t_1(s_1,s_2),\ldots,t_{n-2}(s_1,s_2),s_1,s_2)=0$$
with respect to $s_1$ and setting $(s_1,s_2)=(0,0)$, we have
$$-\frac{\partial z}{\partial s_1}(0) \boldsymbol{\iota^{\prime}}(0)+\sum_{i=1}^{n-2}
\frac{\partial t_i}{\partial s_1}(0) \frac{\partial \boldsymbol{x}}{\partial t_i}(0)
+\frac{\partial \boldsymbol{x}}{\partial t_{n-1}}(0)=0.$$
But by the transversality, $\boldsymbol{\iota^{\prime}}(0)$ is not contained in
${\rm Span}\left( \partial \boldsymbol{x}/\partial t_i(0)\right)_{1\le i\le n}$.
Thus we have $\partial z/\partial s_1(0)=0$ and the identity
\begin{equation}
\label{eq:2-4-2}
\sum_{i=1}^{n-2}\frac{\partial t_i}{\partial s_1}(0)
\frac{\partial \boldsymbol{x}}{\partial t_i}(0)
+\frac{\partial \boldsymbol{x}}{\partial t_{n-1}}(0)=0.
\end{equation}
Similarly, we have $\partial z/\partial s_2(0)=0$ and the identity
\begin{equation}
\label{eq:2-4-3}
\sum_{i=1}^{n-2}\frac{\partial t_i}{\partial s_2}(0)
\frac{\partial \boldsymbol{x}}{\partial t_i}(0)
+\frac{\partial \boldsymbol{x}}{\partial t_{n}}(0)=0.
\end{equation}
This shows that $(0,x^0)$ is a critical point of $f_{\iota}$. We have proved the first part.

To accomplish the proof, we must show that all the critical points are non-degenerate.
We need to compute the Hessian of $z(s_1,s_2)$ at $(s_1,s_2)=(0,0)$ where $z(s_1,s_2)$ is as in (\ref{eq:2-4-1}).
For this purpose we give a system of local equations for the submanifold $\mathcal{Y}_2$, and
we rephrase the fact that the differential
$(p_2|_{\mathcal{D}})_* \colon T_{(x^0,W_0)}\mathcal{D}=T_{(x^0,W_0)}\mathcal{Y}_2\rightarrow T_{W_0}D_X$
is an isomorphism.

Now take the local coordinates (\ref{eq:2-1-3}) of $\mathcal{W}$. Then for $(x,W)$ in this coordinate
neighborhood, $T_xX=T_x(X\cap \mathcal{U}_0)\subset \mathbb{C}^N$ is spanned by the $n$ vectors
$$\alpha_i=\left(\frac{\partial x_1}{\partial t_i},\ldots,
\frac{\partial x_N}{\partial t_i}\right)^t,\ 1\le i\le n,$$
and $T_xW=T_x(W\cap \mathcal{U}_0)$ is spanned by the $k$ vectors
$$\beta_i=(0,\ldots,\stackrel{i}{1},\ldots,0,F_i^{k+1},\ldots,F_i^N)^t,\ 1\le i\le k.$$
$(x,W)\in \mathcal{Y}_2$ if and only if the linear span of these $n+k=N+1$ vectors is $(N-1)$-dimensional.
On the other hand, the origin $(x^0,W_0)$ is in $\mathcal{Y}_2$ and by the assumption (\ref{eq:2-3-4})
$\alpha_1,\ldots,\alpha_{n-2},\beta_1,\ldots,\beta_k$ are linearly independent at $(x^0,W_0)$.
Therefore, the vanishing of the two determinants
$\det(\alpha_1,\ldots,\alpha_{n-2},\alpha_{n-1},\beta_1,\ldots,\beta_k)$
and $\det(\alpha_1,\ldots,\alpha_{n-2},\alpha_n,\beta_1,\ldots,\beta_k)$
gives a system of local equations for $\mathcal{Y}_2$ near $(x^0,W_0)$.
By elementary transformations of matrices, we see that these determinants are equal up to sign to
$\Phi_{n-1}$ and $\Phi_n$ respectively, where
$$\Phi_{\nu}=\det \left( \frac{\partial \boldsymbol{x}}{\partial t_1}
-\sum_{i=1}^k \frac{\partial x_i}{\partial t_1} \boldsymbol{F}_i, \cdots,
\frac{\partial \boldsymbol{x}}{\partial t_{n-2}}
-\sum_{i=1}^k \frac{\partial x_i}{\partial t_{n-2}} \boldsymbol{F}_i,
\frac{\partial \boldsymbol{x}}{\partial t_{\nu}}
-\sum_{i=1}^k \frac{\partial x_i}{\partial t_{\nu}} \boldsymbol{F}_i \right).$$
Here, $\boldsymbol{F}_i=(F_i^{k+1},\ldots,F_i^N)^t$.
Hence $\mathcal{Y}_2$ is locally given by $\Phi_{n-1}=\Phi_{n}=0$.

Now the fact that $(p_2|_{\mathcal{D}})_* \colon T_{(x^0,W_0)}\mathcal{Y}_2\rightarrow T_{W_0}D_X$
is an isomorphism can be rephrased as: the rank of the Jacobian matrix of
$(\Phi_{n-1},\Phi_n,p_2)$ at $(x^0,W_0)$ is equal to $\dim D_X+2=(k+1)(n-1)+1$.
Again by elementary transformations, this is equivalent to the following

\begin{lem}
\label{lem:2-4-1}
Let $(x^0,\mathcal{W}_0)$ and let $\Phi_{n-1}$ and $\Phi_n$ be as in the above.
Then the rank of the $(n+1)\times n$ matrix
\begin{equation}
\label{eq:2-4-4}
\left( \begin{array}{ccc}
\displaystyle\frac{\partial \Phi_{n-1}}{\partial t_1}(0) & \cdots  &
\displaystyle\frac{\partial \Phi_{n-1}}{\partial t_n}(0) \\
\displaystyle\frac{\partial \Phi_{n}}{\partial t_1}(0) & \cdots  &
\displaystyle\frac{\partial \Phi_{n}}{\partial t_n}(0) \\
\displaystyle\frac{\partial \boldsymbol{x}}{\partial t_1}(0) & \cdots &
\displaystyle\frac{\partial \boldsymbol{x}}{\partial t_n}(0) \\
\end{array} \right)
\end{equation}
is equal to $(k+1)(n-1)+1-k(n-1)=n$.
\end{lem}
We perform the following two elementary transformation to (\ref{eq:2-4-4}): let $C_i$ be the $i$-th column
of (\ref{eq:2-4-4}), then 1) add $\sum_{i=1}^{n-2}\partial t_i/\partial s_1(0)C_i$
to the ($n-$1)-th column, and
2) add $\sum_{i=1}^{n-2}\partial t_i/\partial s_2(0)C_i$
to the $n$-th column. Then by (\ref{eq:2-4-2}) and (\ref{eq:2-4-3}), (\ref{eq:2-4-4}) is transformed into
$$\left( \begin{array}{ccccc}
\displaystyle\frac{\partial \Phi_{n-1}}{\partial t_1}(0) & \cdots  &
\displaystyle\frac{\partial \Phi_{n-1}}{\partial t_{n-2}}(0) & A_{11} & A_{12} \\
\displaystyle\frac{\partial \Phi_{n}}{\partial t_1}(0) & \cdots  &
\displaystyle\frac{\partial \Phi_{n}}{\partial t_{n-2}}(0) & A_{21} & A_{22} \\
\displaystyle\frac{\partial \boldsymbol{x}}{\partial t_1}(0) & \cdots &
\displaystyle\frac{\partial \boldsymbol{x}}{\partial t_{n-2}}(0) & \boldsymbol{0} & \boldsymbol{0} \\
\end{array} \right),$$
where
\begin{equation}
\label{eq:2-4-5}
A_{\lambda \mu}=\sum_{i=1}^{n-2}\frac{\partial t_i}{\partial s_{\mu}}(0)
\frac{\partial \Phi_{n+\lambda -2}}{\partial t_i}(0)
+\frac{\partial \Phi_{n+\lambda -2}}{\partial t_{n+\mu -2}}(0).
\end{equation}

Now combining Lemma \ref{lem:2-4-1} and (\ref{eq:2-3-4}), we see that
\begin{equation}
\label{eq:2-4-6}
\det \left( \begin{array}{cc} A_{11} & A_{12} \\
A_{21} & A_{22} \end{array} \right)\neq 0.
\end{equation}

\begin{lem}
\label{lem:2-4-2}
We have the equality
$$\left( \begin{array}{cc} A_{11} & A_{12} \\
A_{21} & A_{22} \end{array} \right)=A_0 \left( \begin{array}{cc}
\displaystyle\frac{\partial^2 z}{{\partial s_1}^2}(0) &
\displaystyle\frac{\partial^2 z}{\partial s_2 \partial s_1}(0) \\
\displaystyle\frac{\partial^2 z}{\partial s_1 \partial s_2}(0) &
\displaystyle\frac{\partial^2 z}{{\partial s_2}^2}(0) \\
\end{array} \right),$$
where $A_0=\det\left(\partial \boldsymbol{x}/\partial t_1(0),\ldots,
\partial \boldsymbol{x}/\partial t_{n-2}(0),\boldsymbol{\iota^{\prime}}(0) \right)$.
\end{lem}

By (\ref{eq:2-3-8}), (\ref{eq:2-4-6}), and Lemma \ref{lem:2-4-2}, it follows that the Hessian of $z(s_1,s_2)$
at $(s_1,s_2)=(0,0)$ is non-zero. Thus, $(0,x^0)$ is a non-degenerate critical point of $f_{\iota}$.
This completes the proof of Theorem \ref{thm:2-3-4}, modulo Lemma \ref{lem:2-4-2}.

Lemma \ref{lem:2-4-2} can be proved by a straightforward computation.
We only give an outline of the proof of
$A_{11}=A_0 \partial^2 z/{\partial s_1}^2(0)$.
The first claim is that
\begin{equation}
\label{eq:2-4-7}
\frac{\partial \Phi_{n-1}}{\partial t_i}(0)=
\det \left( \frac{\partial \boldsymbol{x}}{\partial t_1}(0),\ldots,
\frac{\partial \boldsymbol{x}}{\partial t_{n-2}}(0),
\sum_{\ell=1}^{n-2}\frac{\partial t_{\ell}}{\partial s_1}(0)
\frac{\partial^2 \boldsymbol{x}}{\partial t_i \partial t_{\ell}}(0)+
\frac{\partial^2 \boldsymbol{x}}{\partial t_i \partial t_{n-1}}(0) \right),
\end{equation}
which can be proved by using (\ref{eq:2-4-2}).

For $k+1\le j\le N$, differentiating twice the identity
$$\psi_j(z(s_1,s_2),t_1(s_1,s_2),\ldots,t_{n-2}(s_1,s_2),s_1,s_2)=0$$
with respect to $s_1$ and setting $(s_1,s_2)=(0,0)$, we have
\begin{equation}
\label{eq:2-4-8}
-\frac{d\iota_0^j}{dz}(0)\frac{\partial^2 z}{{\partial s_1}^2}(0)+
\sum_{i=1}^n \left( \left(\sum_{h=1}^n \frac{\partial^2 x_j}{\partial t_h \partial t_i}(0)
\frac{\partial t_h}{\partial s_1}(0)\right)\frac{\partial t_i}{\partial s_1}(0)+
\frac{\partial x_j}{\partial t_i}(0)\frac{\partial^2 t_i}{{\partial s_1}^2}(0) \right)=0
\end{equation}
(note that $\partial z/\partial s_1(0)=\partial z/\partial s_2(0)=0$,
$\partial t_{n-1}/\partial s_1(0)=1$, and $\partial t_n/\partial s_1(0)=0$).

Using (\ref{eq:2-4-5}), (\ref{eq:2-4-7}), and (\ref{eq:2-4-8}),
we can get the desired formula.

\section{The Meyer function for a projective variety}

\subsection{Main theorem}
Recall the situation arising from $X\subset \mathbb{P}_N$ as described in section 1.
We focus on the topological monodromy
$$\rho_X\colon \pi_1(U^X)\rightarrow \Gamma_g$$
of (\ref{eq:1-1}).

\begin{thm}
\label{thm:3-1-1}
There exists a uniquely determined $\mathbb{Q}$-valued 1-cochain
$\phi_X \colon \pi_1(U^X)\rightarrow \mathbb{Q}$ whose coboundary equals to
$\rho_X^*\tau_g$. In particular, we have $\rho_X^*[\tau_g]=0\in H^2(\pi_1(U^X);\mathbb{Q})$.
\end{thm}

Here we comment about the group $\pi_1(U^X)$. If the codimension of $D_X$ is $\ge 2$, $\pi_1(U^X)$ is trivial
since $G_k(\mathbb{P}_N)$ is simply connected. Suppose the codimension of $D_X$ is 1.
Then $\pi_1(U^X)$ is finitely presentable since $U^X$ is an affine algebraic variety,
and the first Betti number $b_1(\pi_1(U^X))$ is zero (see Lemma \ref{lem:3-3-1}).
Moreover, $\pi_1(U^X)$ is normally generated by a single element, called a \textit{lasso}.
Roughly speaking, a lasso is an element of $\pi_1(U^X)$ going once around $D_X$.
The precise definition is as follows. We fix some base point in $U^X$.
Let $W_0$ be a smooth point of $D_X$ and $(z_1,\ldots,z_m)$
be local coordinates of $G_k(\mathbb{P}_N)$ centered at $W_0$, such that $D_X$ is locally given by
$z_1=0$. For a sufficiently small $\varepsilon >0$, consider the loop
$$[0,1]\rightarrow U^X,\ t\mapsto (\varepsilon e^{2\pi \sqrt{-1}t},0,\ldots,0)$$
defined in this coordinate neighborhood. Joining this loop with a path in $U^X$ from the base point of $U^X$
to $(\varepsilon,0,\ldots,0)$, we get an element of $\pi_1(U^X)$, which is called a lasso around $D_X$.
The irreducibility of $D_X$ implies that all lassos are conjugate to each other and $\pi_1(G_k(\mathbb{P}_N))=1$
implies that $\pi_1(U^X)$ is normally generated by a lasso. Since $\phi_X$ is a class function (see section 1),
the value of $\phi_X$ on any lasso is constant. This value can be computed from various invariants of $X$.
For details, see subsection 3.5.
Nevertheless, the values of $\phi_X$ on an element other than lasso seems difficult to know.

The proof of Theorem \ref{thm:3-1-1} will be given in the next two subsections.
The following argument is a generalization of sections 3 and 4 of \cite{Ku}.

\subsection{Proof of the existence}
It suffices to consider the case when the codimension of $D_X$ is 1.
The existence of $\phi_X$ is equivalent to $\rho_X^*[\tau_g]=0\in H^2(\pi_1(U^X);\mathbb{Q})$.

We first embed $H^2(\pi_1(U^X);\mathbb{Q})$ into another space. Let $V_{k+1,N+1}$ be the (complex) Stiefel manifold
of all $(k+1)$-frames of $\mathbb{C}^{N+1}$. Regarding $\mathbb{P}_N$ as the projectivization of $\mathbb{C}^{N+1}$
we have the natural projection $q \colon V_{k+1,N+1}\rightarrow G_k(\mathbb{P}_N)$, which is a principal
$GL(k+1;\mathbb{C})$ bundle. Let $\widetilde{D}_X=q^{-1}(D_X)$,
$\widetilde{U}^X=V_{k+1,N+1} \setminus \widetilde{D}_X$, and $\widetilde{E}=q^{-1}(E)$.
For simplicity we use the same letter $q$ for the restriction
$q|_{\widetilde{U}^X}\colon \widetilde{U}^X\rightarrow U^X$.

We have the following short exact sequence with rational coefficients:
$$H^0(U^X)\stackrel{\cup c_1}{\longrightarrow}H^2(U^X)
\stackrel{q^*}{\longrightarrow}H^2(\widetilde{U}^X).$$
This is derived from the 5-term exact sequence of the principal bundle
$q\colon \widetilde{U}^X \rightarrow U^X$. Here, $c_1$ is the first Chern class, which is the restriction of
a generator of $H^2(G_k(\mathbb{P}_N))$ to $H^2(U^X)$. But since $D_X$ is of codimension 1, the first Chern class
$c_1([D_X])\in H^2(G_k(\mathbb{P}_N))$ is defined and is also a generator of $H^2(G_k(\mathbb{P}_N))$.
The point here is that $H^2(G_k(\mathbb{P}_N))$ is of rank 1. Clearly the restriction of $c_1([D_X])$ to
$U^X=G_k(\mathbb{P}_N)\setminus D_X$ is zero, therefore $c_1$ is also zero. Thus, we have the injective homomorphism
\begin{equation}
\label{eq:3-2-1}
q^*\colon H^2(U^X;\mathbb{Q})\hookrightarrow H^2(\widetilde{U}^X;\mathbb{Q}).
\end{equation}
Let $\chi\colon \pi_1(\widetilde{U}^X)\rightarrow \pi_1(U^X)$ be the homomorphism between fundamental groups
induced by $q$. Since for any space $X$ there is the natural injection $H^2(\pi_1(X))\rightarrow H^2(X)$
of the second cohomology with arbitrary coefficients, (\ref{eq:3-2-1}) implies
that we have the injective homomorphism
\begin{equation}
\label{eq:3-2-2}
\chi^*\colon H^2(\pi_1(U^X);\mathbb{Q})\hookrightarrow H^2(\pi_1(\widetilde{U}^X);\mathbb{Q})
\end{equation}
induced by $\chi$.

Next we show $\chi^*\rho_X^*[\tau_g]=0\in H^2(\pi_1(\widetilde{U}^X);\mathbb{Z})$. Let
$$\widetilde{\mathcal{W}}:=\left\{ (x,\widetilde{W})\in \mathbb{P}_N \times V_{k+1,N+1};
x\in X\cap q(\widetilde{W})\right\}$$
and $\tilde{p}_2\colon \widetilde{\mathcal{W}}\rightarrow V_{k+1,N+1}$ be the second projection, and
$$\widetilde{\mathcal{C}^X}:=\left\{ (x,\widetilde{W})\in \mathbb{P}_N \times \widetilde{U}^X;
x\in X\cap q(\widetilde{W})\right\}.$$
The second projection $\tilde{p}_X\colon \widetilde{\mathcal{C}^X}\rightarrow \widetilde{U}^X$
is a family of Riemann surfaces, which is the pull back of $p_X\colon \mathcal{C}^X\rightarrow U^X$
by $q$. The associated topological monodromy is $\tilde{\rho}_X:=\rho_X \circ \chi$.

We construct a 1-cochain $c\colon \pi_1(\widetilde{U}^X)\rightarrow \mathbb{Z}$ whose coboundary $\delta c$
coincides with ${\tilde{\rho}_X}^*\tau_g$. The point here is $V_{k+1,N+1}\setminus \widetilde{E}$ is 2-connected.
This follows from the two facts: 1) the Stiefel manifold $V_{k+1,N+1}$ is $2(N-k)$-connected and $2(N-k)=2n-2\ge 2$,
and 2) the complex codimension of $\widetilde{E}\subset V_{k+1,N+1}$ is $\ge 2$ (see Lemma \ref{lem:2-3-1}).
All of the spaces that we consider in the rest of this subsection as well as all of the maps are based,
otherwise stated.

\noindent \textbf{Construction of $c$.}
Let $\ell\colon S^1\rightarrow \widetilde{U}^X$ be a $C^{\infty}$-loop, i.e., a $C^{\infty}$-map from
$S^1$ to $\widetilde{U}^X$. Since $V_{k+1,N+1}\setminus \widetilde{E}$ is simply connected we can extend
$\ell$ to a $C^{\infty}$-map $\tilde{\ell}\colon D^2\rightarrow V_{k+1,N+1}\setminus \widetilde{E}$
which is transverse to $\widetilde{D}_X$. Here we make the identifications $S^1=\{z \in \mathbb{C};|z|=1 \}$
and $D^2=\{z\in \mathbb{C};|z|\le 1\}$, and endow them the usual orientation: the orientation of $D^2$
is induced by that of $\mathbb{C}$ and $S^1$ goes around $D^2$ by counter clockwise manner.

By Proposition \ref{prop:2-3-3} the pull back
$\tilde{\ell}^*\widetilde{\mathcal{W}}:=(q \circ \tilde{\ell})^*\mathcal{W}$ has the natural structure
of a compact oriented 4-dimensional $C^{\infty}$-manifold with boundary. The orientation is induced
by the orientation of $D^2$ and that of the general fibers of $\tilde{\ell}^*\widetilde{\mathcal{W}}$,
which have the natural orientations as compact Riemann surfaces. Set
$$c([\ell]):={\rm Sign}(\tilde{\ell}^*\widetilde{\mathcal{W}})\in \mathbb{Z}.$$
Here $[\ell]\in \pi_1(\widetilde{U}^X)$ is the element represented by $\ell$,
and the right hand side is the signature of $\tilde{\ell}^*\widetilde{\mathcal{W}}$.

\begin{prop}
\label{prop:3-2-1}
The above definition of $c$ is well defined. The 1-cochain $c$ is a class function on $\pi_1(\widetilde{U}^X)$
and $c(x^{-1})=-c(x)$ for $x\in \pi_1(\widetilde{U}^X)$. We have $\delta c=-{\tilde{\rho}_X}^*\tau_g$.
\end{prop}

\begin{proof}
Let $\ell_0$ and $\ell_1$ be $C^{\infty}$-loops in $\widetilde{U}^X$. Suppose that the elements of
$\pi_1(\widetilde{U}^X)$ represented by them are conjugate to each other. Then there exists a
$C^{\infty}$-homotopy $H\colon S^1 \times [0,1]\rightarrow \widetilde{U}^X$ such that $H(\cdot,0)=\ell_0$
and $H(\cdot,1)=\ell_1$ (caution: we do not require that $H(\cdot,t)$ is a base preserving map
for every $t\in [0,1]$). Identify the 2-sphere $S^2$ as
$$S^2\cong (S^1\times [0,1])\cup (D^2\times \{0\}) \cup (D^2\times \{1\})$$
and take some extensions $\tilde{\ell}_i\colon D^2\times \{i\}\rightarrow V_{k+1,N+1}\setminus \widetilde{E}$
of $\ell_i$ for $i=0,1$. Then piecing $H$, $\tilde{\ell}_0$, and $\tilde{\ell}_1$ together,
we can construct a $C^{\infty}$-map $\widetilde{H}\colon S^2\rightarrow V_{k+1,N+1}\setminus \widetilde{E}$
which is transverse to $\widetilde{D}_X$. Introduce the orientation of $S^2$ such that
$D^2\times \{0\} \hookrightarrow S^2$ is orientation preserving. Then $D^2\times \{1\}\hookrightarrow S^2$
is orientation reversing and the pull back $\widetilde{H}^*\widetilde{\mathcal{W}}$ is a closed
oriented 4-dimensional $C^{\infty}$-manifold. Moreover, since $\pi_2(V_{k+1,N+1}\setminus \widetilde{E})=0$,
$\widetilde{H}$ extends to a $C^{\infty}$-map from the 3-ball to $V_{k+1,N+1}\setminus \widetilde{E}$
which is transverse to $\widetilde{D}_X$. Hence $\widetilde{H}^*\widetilde{\mathcal{W}}$
is the boundary of a 5-dimensional manifold and the signature of $\widetilde{H}^*\widetilde{\mathcal{W}}$
is zero. Now by the Novikov additivity of the signature we have
$$0={\rm Sign}(\widetilde{H}^*\widetilde{\mathcal{W}})={\rm Sign}(\tilde{\ell}_0^*\widetilde{\mathcal{W}})
-{\rm Sign}(\tilde{\ell}_1^*\widetilde{\mathcal{W}}).$$
This proves that $c$ is well defined and $c$ is a class function, i.e., $c(xyx^{-1})=c(y)$
for $x,y\in \pi_1(\widetilde{U}^X)$. Since changing the orientation of a manifold changes
the sign of its signature, the property $c(x^{-1})=-c(x)$ is clear.

We next prove that $\delta c=-{\tilde{\rho}_X}^*\tau_g$, i.e.,
\begin{equation}
\label{eq:3-2-3}
c([\ell_0])+c([\ell_1])-c([\ell_0][\ell_1])=-{\tilde{\rho}_X}^*\tau_g([\ell_0],[\ell_1])
\end{equation}
for any based $C^{\infty}$-loops $\ell_0$ and $\ell_1$. Let $D_i$, $0\le i\le 2$, be embedded three disjoint
closed 2-disks in $S^2$ and we denote its boundary circle by $S^1_i$. Let
$P:=S^2 \setminus \coprod_{i=0}^2 {\rm Int}(D_i)$. Since $P$ has the homotopy type of the bouquet $S^1 \vee S^1$,
we can construct a $C^{\infty}$-map $L\colon P\rightarrow \widetilde{U}^X$ such that the restriction of
$L$ to $S^1_i\cong S^1$ is equal to $\ell_i$ for $i=0,1$. Then the restriction of $L$ to $S^1_2$ is homotopic to
the inverse of the composition loop $\ell_0\cdot \ell_1$. Notice that ${\rm Sign}(L^*\widetilde{\mathcal{W}})$
is equal to ${\tilde{\rho}_X}^*\tau_g([\ell_0],[\ell_1])$. Take some extensions $\tilde{\ell}_0$, $\tilde{\ell}_1$,
and $\widetilde{\ell_0\cdot \ell_1}$ of $\ell_0$, $\ell_1$, and $\ell_0\cdot \ell_1$, respectively.
Then by piecing them and $L$ together we have a $C^{\infty}$-map
$\widetilde{L}\colon S^2\rightarrow V_{k+1,N+1}\setminus \widetilde{E}$
which is transverse to $D_X$. Again, the vanishing of $\pi_2(V_{k+1,N+1}\setminus \widetilde{E})$
implies that the signature of $\widetilde{L}^*\widetilde{\mathcal{W}}$ is zero.
Finally, by the Novikov additivity we have
$$0={\rm Sign}(\widetilde{L}^*\widetilde{\mathcal{W}})=
{\rm Sign}(\tilde{\ell}_0^*\widetilde{\mathcal{W}})+{\rm Sign}(\tilde{\ell}_0^*\widetilde{\mathcal{W}})
-{\rm Sign}(\widetilde{\ell_0\cdot \ell_1}^*\widetilde{\mathcal{W}})+{\rm Sign}(L^*\widetilde{\mathcal{W}}),$$
but this equation is equivalent to (\ref{eq:3-2-3}). This completes the proof. 
\end{proof}

By Proposition \ref{prop:3-2-1}, we have $\chi^*\rho_X^*[\tau_g]=0\in H^2(\pi_1(\widetilde{U}^X);\mathbb{Z})$.
Combining this with the injection (\ref{eq:3-2-2}) we have $\rho_X^*[\tau_g]=0\in H^2(\pi_1(U^X);\mathbb{Q})$.
This completes the proof of the existence of $\phi_X$.

\subsection{Proof of the uniqueness}
The uniqueness of $\phi_X$ follows from the following

\begin{lem}
\label{lem:3-3-1}
The first cohomology group of $\pi_1(U^X)$ is trivial over rationals: 
$$H^1(\pi_1(U^X);\mathbb{Q})={\rm Hom}(\pi_1(U^X),\mathbb{Q})=0.$$
\end{lem}

\begin{proof}
It suffices to consider the case when the codimension of $D_X$ is 1.
Consider the following commutative diagram among (co)homology groups with integer coefficients:
$$\xymatrix{
\mathbb{Z} \cong H_2(G_k(\mathbb{P}_N)) \ar[d]_{\cong} \ar[r] &
H_2(G_k(\mathbb{P}_N),U^X) \ar[d]^{\cong} \ar[r]
& H_1(U^X) \ar[r] & 0 \\
H^{2\dim G_k(\mathbb{P}_N)-2}(G_k(\mathbb{P}_N)) \ar[r]^{j^*} &
H^{2\dim G_k(\mathbb{P}_N)-2}(D_X)\cong \mathbb{Z}. \\}$$
The vertical isomorphisms are Poincar\'e duality. Note that $H^{2\dim G_k(\mathbb{P}_N)-2}(D_X)\cong \mathbb{Z}$
since $D_X$ is irreducible. The first horizontal sequence is exact and is a part of the homology sequence of
the pair $(G_k(\mathbb{P}_N),U^X)$ and $j^*$ is induced by the inclusion $D_X \hookrightarrow G_k(\mathbb{P}_N)$.
Then the generator of $H_2(G_k(\mathbb{P}_N))$ is mapped to a positive integer times the generator of
$H^{2\dim G_k(\mathbb{P}_N)-2}(D_X)$, the fundamental class of $D_X$ (this positive integer is denoted
by $\deg D_X$ and will be studied in the next subsection). Thus $H_1(U^X)$, which is isomorphic to
the abelianization of $\pi_1(U^X)$, is a cyclic group of finite order. This completes the proof.
\end{proof}

Now Theorem \ref{thm:3-1-1} is established.

\subsection{Theory of Lefschetz pencils}
In this subsection we recall the definition of the degree of an analytic subset in a Grassmannian
and describe a method to compute the degree of $D_X\subset G_k(\mathbb{P}_N)$.
This will be used to compute the value of $\phi_X$ on a lasso.

First we treat the case of classical dual varieties, namely when $n=2$. Let $X\subset \mathbb{P}_N$ be a smooth
projective surface. Then $G_k(\mathbb{P}_N)=G_{N-1}(\mathbb{P}_N)$ is nothing but the dual projective space
$\mathbb{P}_N^{\vee}$ and $D_X$ is the dual variety of $X$. Let $L$ be a line of $\mathbb{P}_N^{\vee}$ avoiding
the singular points of $D_X$ and meeting $D_X$ transversally. Note that generic lines of $\mathbb{P}_N^{\vee}$
satisfy this condition. We denote by $i_L$ the inclusion $L\hookrightarrow \mathbb{P}_N^{\vee}$.
Then as explained in \cite{Katz} or \cite{La}, especially (1.6.3) of \cite{La}, the projection
$i_L^*\mathcal{W}\rightarrow L$ is a \textit{holomorphic Lefschetz fibration} with the set of
critical values being $L\cap D_X$, in the following sense.

\begin{dfn}
\label{dfn:3-4-1}
Let $Y$ be a complex surface and $C$ a compact Riemann surface. A proper surjective holomorphic map
$f\colon Y\rightarrow C$ is called a holomorphic Lefschetz fibration if the number of critical values of $f$
is finite and over each critical value, there exists only one critical point near which
$f$ locally looks like $(z_1,z_2)\mapsto z_1^2+z_2^2$.
\end{dfn}
In particular $\deg D_X$, which is equal to the number $\sharp(L\cap D_X)$ by definition,
is equal to the number of critical points of $i_L^*\mathcal{W}\rightarrow L$.

More generally for a projective variety $X\subset \mathbb{P}_N$, we can take a generic line $L$ of 
$\mathbb{P}_N^{\vee}$ and consider the family of hyperplane sections $\{ H\cap X\}_{H\in L}$ of $X$, parametrized
by $L$. This construction is called \textit{Lefschetz pencils} and very useful to study the topology of $X$.

We slightly generalize the above construction to the case of general $n$. We shall start by giving the definition
of the \textit{degree of an analytic subset in a Grassmannian} following \cite{GKZ} Chapter 3, section 2-A.

Recall that $G_{p,q}$ is the Grassmannian of $p$-planes of $\mathbb{C}^q$.
By a \textit{line of $G_{p,q}$} is meant a curve embedded
in $G_{p,q}$ which can be written as
\begin{equation}
\label{eq:3-4-1}
P_{NM}:=\{ W\in G_{p,q};N\subset W\subset M \}
\end{equation}
for some $(p-1)$-plane $N$ and $(p+1)$-plane $M$ satisfying $N\subset M$.

Let $Z$ be an analytic subset of $G_{p,q}$. Take a line $P_{MN}$ of $G_{p,q}$ such that
\begin{enumerate}
\item $P_{MN}\cap S(Z)=\emptyset$, where $S(Z)$ denotes the set of singular points of $Z$,
\item $P_{MN}$ and $Z$ meet transversally.
\end{enumerate}
Then the intersection $P_{MN}\cap Z$ consists of finitely many points. We define $\deg Z \in \mathbb{Z}$ by
$$\deg Z:=\sharp (P_{MN}\cap Z)$$
where $P_{MN}$ satisfies the above two conditions. When $Z$ is a hypersurface this number is positive
(see \cite{GH} p.64), and has the following topological interpretation. Let $c_1([Z])\in H^2(G_{p,q};\mathbb{Z})$
be the first Chern class of the line bundle over $G_{p,q}$ determined by $Z$ and $[P_{NM}]\in H_2(G_{p,q};\mathbb{Z})$
the homology class represented by the embedded 1-dimensional projective space $P_{NM}$. Note that $[P_{NM}]$
is a generator of $H_2(G_{p,q};\mathbb{Z})\cong\mathbb{Z}$. Then $\deg Z$ is equal to the result of the Kronecker
pairing $\langle c_1([Z]),[P_{NM}]\rangle$. If $n=2$, this degree coincides with the usual definition of the
degree of a projective hypersurface in $G_{N-1}(\mathbb{P}_N)=\mathbb{P}_N^{\vee}$.
When the codimension of $Z$ is $\ge 2$, $\deg Z=0$.

We remark that generic lines of $G_{p,q}$ satisfy the above two conditions in the following sense.
Let us consider the space parametrizing all lines of $G_{p,q}$; namely, let
$$\mathcal{L}_{p,q}:=\left\{ (N,M)\in G_{p-1,q}\times G_{p+1,q}; N\subset M \right\}.$$
Then the set of $(N,M)\in \mathcal{L}_{p,q}$ such that $P_{NM}$ satisfying the above two conditions is
non-empty and Zariski open in $\mathcal{L}_{p,q}$. This is proved by an application of Sard's lemma for
varieties to the second projection
$\{ (z,(N,M))\in Z\times \mathcal{L}_{p,q}; z\in P_{NM}\}\rightarrow \mathcal{L}_{p,q}$.

Let us return to our setting: $X\subset \mathbb{P}_N$ is a $n$-dimensional smooth projective variety,
$D_X$ is the $k$-th associated variety. Let $E_X=E\subset D_X$ be as in Lemma \ref{lem:2-3-1} when
the codimension of $D_X$ is 1, $E_X=\emptyset$ when the codimension of $D_X$ is $\ge 2$.
Here we put the subscript $_X$ to $E$ to indicate its dependence on $X$. Let
$$\mathcal{L}_k(\mathbb{P}_N):=
\left\{ (N,M)\in G_{k-1}(\mathbb{P}_N)\times G_{k+1}(\mathbb{P}_N); N\subset M\right\};$$
this is clearly isomorphic to $\mathcal{L}_{k+1,N+1}$. For $(N,M)\in \mathcal{L}_k(\mathbb{P}_N)$,
the corresponding line $P_{NM}\subset G_k(\mathbb{P}_N)$ is defined by the same way as (\ref{eq:3-4-1}).

The existence of Lefschetz pencils for the case of general $n$ is stated as follows.

\begin{thm}[Existence of Lefschetz pencils]
\label{thm:3-4-2}
Let $P_{NM}$ be a line of $G_k(\mathbb{P}_N)$ not meeting $E_X$ and meeting $D_X$ transversally. Then the projection
$i_{NM}^*\mathcal{W}\rightarrow P_{NM}$, where $i_{NM}\colon P_{NM}\hookrightarrow G_k(\mathbb{P}_N)$ denotes
the inclusion, is a holomorphic Lefschetz fibration in the sense of Definition \ref{dfn:3-4-1}. Moreover
$\deg D_X$ is equal to the number of critical points of $i_{NM}^*\mathcal{W}\rightarrow P_{NM}$.
\end{thm}

By Theorem \ref{thm:2-3-4}, the remaining to show is the number of critical points over each critical value is
just one. If $n=2$, there is nothing to prove as remarked before Definition \ref{dfn:3-4-1}.
We only remark that in the proof of (1.6.3) of \cite{La}, the bi-duality theorem plays a key role.

To reduce the case of general $n$ to the case of $n=2$, we will cut $X$ with a generic $(k+1)$-plane. The result
will be a smooth projective surface in the $(k+1)$-plane. We prepare some notations. $P_{NM}$ can be considered
as a line of $M^{\vee}$, the dual projective space of $M$. Then we write it by $L_N$. For $M\in G_k(\mathbb{P}_N)$,
let $X^{\prime}:=M\cap X\subset M$. If $M$ meets $X$ transversally, $X^{\prime}$ is a smooth surface in $M$.
Then we can consider $D_{X^{\prime}}$ and $E_{X^{\prime}}$ in $M^{\vee}$.

\begin{lem}
\label{lem:3-4-3}
There exists a point $(N,M)\in \mathcal{L}_k(\mathbb{P}_N)$ such that:
\begin{enumerate}
\item the $(k+1)$-plane $M$ meets $X$ transversally (hence $X^{\prime}$ is a smooth projective surface).
\item the line $P_{NM}$ does not meet $E_X$ and meets $D_X$ transversally.
\item the line $L_N$ does not meet $E_{X^{\prime}}$ and meets $D_{X^{\prime}}$ transversally.
\end{enumerate}
\end{lem}
\begin{proof}
The set of points in $\mathcal{L}_k(\mathbb{P}_N)$ satisfying the conditions 1 and 2 is non-empty and Zariski open
in $\mathcal{L}_k(\mathbb{P}_N)$. Let $(N^{\prime},M)$ be a point in this set. Since the set of lines of $M^{\vee}$
not meeting $E_{X^{\prime}}$ and meeting $D_{X^{\prime}}$ transversally is non-empty and Zariski open in the space
of all lines of $M^{\vee}$, there exists a line $L_N$ near $L_{N^{\prime}}$ such that $(N,M)$ satisfies all the
three conditions.
\end{proof}

\noindent \textbf{Implications of the Lemma.} Let $(N,M)$ be as in Lemma \ref{lem:3-4-3}. We have the natural
inclusion $\iota_M\colon M^{\vee}\hookrightarrow G_k(\mathbb{P}_N)$. Since $M$ meets $X$ transversally, for
$H\in M^{\vee}$ the conditions $H\in D_{X^{\prime}}$ and $\iota_M(H)\in D_X$ are equivalent.
Therefore we have the injection
\begin{equation}
\label{eq:3-4-2}
\iota_M|_{U^{X^{\prime}}}\colon U^{X^{\prime}}\hookrightarrow U^X,
\end{equation}
where $U^{X^{\prime}}=M^{\vee}\setminus D^{X^{\prime}}$, and the bijection
$\iota_M|_{L_N \cap D_{X^{\prime}}}\colon L_N \cap D_{X^{\prime}}\stackrel{\cong}{\rightarrow} P_{NM}\cap D_X$.
In particular, $D_X$ is a hypersurface if and only if $D_{X^{\prime}}$ is a hypersurface and
we have $\deg D_X=\deg D_{X^{\prime}}$.

For simplicity we identify $L_N \cap D_{X^{\prime}}$ with $P_{NM}\cap D_X$ and write it by $D_{NM}$.
Let $U_{NM}=P_{NM} \setminus D_{NM}$. By the inclusion $U_{NM} \hookrightarrow U^X$
(resp. $U_{NM} \hookrightarrow U^{X^{\prime}}$), any loop in $U_{NM}$ going once around a point of $D_{NM}$
is mapped to a lasso around $D_X$ (resp. $D_{X^{\prime}}$), hence a lasso around $D_{X^{\prime}}$ is mapped
to a lasso around $D_X$ by the map (\ref{eq:3-4-2}). Consider the group homomorphism
$j_M \colon \pi_1(U^{X^{\prime}})\rightarrow \pi_1(U^X)$ induced by (\ref{eq:3-4-2}). The uniqueness of
$\phi_{X^{\prime}}$ shows that $j_M^*\phi_X$ coincides with $\phi_{X^{\prime}}$. Thus the value of $\phi_X$ on
a lasso around $D_X$ coincides with the value of $\phi_{X^{\prime}}$ on a lasso around $D_{X^{\prime}}$.

In this way we can reduce the computation of $\deg D_X$ or the value of $\phi_X$ on a lasso around $D_X$
to the case of $n=2$.

\begin{proof}[Proof of Theorem \ref{thm:3-4-2}]
Let $(N,M)$ be as in Lemma \ref{lem:3-4-3}.
Let $i_N^{\prime}\colon L_N\rightarrow M^{\vee}$ be the inclusion and
$\mathcal{W}^{\prime}:=\{ (x,W)\in M \times M^{\vee};x\in X^{\prime}\cap W\}$.
We can consider the pull back $(i_N^{\prime})^*{\mathcal{W}^{\prime}}$. Since $\dim X^{\prime}=2$ the remark right
after the statement of Theorem \ref{thm:3-4-2} applies, so $(i_N^{\prime})^*{\mathcal{W}^{\prime}}\rightarrow L_N$
is a holomorphic Lefschetz fibration. Therefore, $i_{NM}^*\mathcal{W}\rightarrow P_{NM}$ is also a holomorphic
Lefschetz fibration because $\iota_M$ induces the isomorphism
$$(i_N^{\prime})^*{\mathcal{W}^{\prime}}\stackrel{\cong}{\rightarrow} i_{NM}^*\mathcal{W}$$
between the families of algebraic curves over $L_N=P_{NM}$. Thus we have proved that:
there exists a line $P_0=P_{NM}$ of $G_k(\mathbb{P}_N)$ not meeting $E_X$ and meeting $D_X$ transversally
such that the projection $i_{NM}^*\mathcal{W}\rightarrow P_{NM}$ is a holomorphic Lefschetz fibration.

Let $P_1$ be a line not meeting $E_X$ and meeting $D_X$ transversally. For $j=0,1$, we denote the inclusion map
$P_j\hookrightarrow G_k(\mathbb{P}_N)$ by $i_j$. Let $\mathcal{V}$ be the space of lines of $G_k(\mathbb{P}_N)$
not meeting $E_X$ and meeting $D_X$ transversally. This is non-empty and Zariski open, hence connected. Thus there
exists a differentiable path in $\mathcal{V}$ joining $P_0$ and $P_1$, inducing a deformation equivalence of class
$C^{\infty}$ between $i_0^*\mathcal{W}$ and $i_1^*\mathcal{W}$ as a family of algebraic curves over 1-dimensional
projective space. We already know  $i_0^*\mathcal{W}\rightarrow P_0$ is a holomorphic Lefschetz fibration, so
$i_1^*\mathcal{W}\rightarrow P_1$ is also a holomorphic Lefschetz fibration.
\end{proof}

Taking into account that $D_X\setminus E_X$ is connected, similar argument shows the following

\begin{cor}
\label{cor:3-4-4}
In Theorem \ref{thm:2-3-4}, the number of critical points of $f_{\iota}\colon \iota^*\mathcal{W}\rightarrow \Delta$
is 1. The singular fiber $f_{\iota}^{-1}(0)$ has just one nodal singularity and its topological type does not depend
on the choice of $\iota$.
\end{cor}

\subsection{Computations}
In this subsection we will give a formula for the value of $\phi_X$ on a lasso around $D_X$ from the data of
various invariants of $X$. In view of the discussion 'implication of the lemma' in subsection 3.4, we may focus
on the case $n=2$.

First we review theory of Lefschetz pencils following \cite{La}.
Let $X\subset \mathbb{P}_N$ be a smooth projective surface and $L$ a generic line of $\mathbb{P}_N^{\vee}$,
as in the beginning of subsection 3.4. Let $D_L=L\cap D_X$, $U_L=L \setminus D_L$ and $X_H=H\cap X$ for $H\in L$.
We have a family of Riemann surfaces over $U_L$ by restricting $i_L^*\mathcal{W}\rightarrow L$. Choosing a base point
$H_0\in U_L$, let $\rho_X$ be the associated topological monodromy. The source of $\rho_X$ is
the fundamental group $\pi_1(U_L)=\pi_1(U_L,H_0)$.

\noindent \textbf{Picard-Lefschetz formula.}
For $H^{\prime}\in D_L$ choose a path $\ell$ from $H_0$ to $H^{\prime}$ and let $\sigma$ be the element of $\pi_1(U_L)$
represented by a loop going to a point nearby $H^{\prime}$ along $\ell$, then going once around $H^{\prime}$ by
counter clockwise manner and then coming back along $\ell$. Then the famous Picard-Lefschetz formula says that
$\rho_X(\sigma)$ is the \textit{inverse} (recall our conventions about monodromies) of
the right hand Dehn twist along some simple closed curve $C_{\sigma}$,
called the \textit{vanishing cycle}, on $X_{H_0}$. The adjective 'vanishing' comes in because looking
at the fiber $H\cap X$ when $H$ moves along $\ell$, $X_{H^{\prime}}$ looks like obtained from $X_{H_0}$
by pinching $C_{\sigma}$ into a point.

Let $V$ be the submodule of $H_1(X_{H_0})$ generated by all the vanishing cycles. Homology with coefficients
in some principal ideal domain is considered. Then the equality
\begin{equation}
\label{eq:3-5-1}
V={\rm Ker}(i_*\colon H_1(X_{H_0})\rightarrow H_1(X)),
\end{equation}
where $i_*$ is induced by the inclusion $X_{H_0}\hookrightarrow X$, holds. See \cite{La}, (3.8.2). In particular, if
$H_1(X)=0$, the vanishing cycles generate the homology of the reference
fiber $X_{H_0}$.

\begin{lem}
\label{lem:3-5-1}

Suppose the genus of $X_{H_0}$ is positive and for some principal ideal domain $R$, the rank
of $H_1(X)=H_1(X;R)$ is less than twice the genus of $X_{H_0}$. Then every
singular fiber of $i_L^*\mathcal{W}\rightarrow L$, i.e., the inverse image of of a point of $D_L$, is irreducible.
\end{lem}
\begin{proof}
First remark that for any choice of $H^{\prime}$ and $\ell$, $\sigma \in \pi_1(U_L)$ is mapped to a lasso
around $D_X$ by the homomorphism $\pi_1(U_L)\rightarrow \pi_1(U^X)$ induced by the inclusion.
Thus for any two vanishing cycles the Dehn twists along them are conjugate to each other in the mapping
class group of $X_{H_0}$.

Suppose there exists a reducible fiber. This means that there exists a vanishing cycle which is a separating
simple closed curve. Then all the vanishing cycles are separating by the remark above. Since any separating
simple closed curve is zero as a homology class, this implies $V=0$. But by the assumption and
(\ref{eq:3-5-1}) we also have $V=H_1(X_{H_0})\neq 0$, a contradiction. This completes the proof.
\end{proof}

\begin{prop}
\label{prop:3-5-2}
Let $X\subset \mathbb{P}_N$ be a smooth projective surface and
$g$ the genus of a generic hyperplane section $H\cap X$, $H\in U^X$.
Assume that $g>0$ and the rank of $H_1(X;R)$ is less than $2g$ for some principal ideal domain $R$. Suppose $D_X$ is
a hypersurface and let $\sigma_X\in \pi_1(U^X)$ be a lasso around $D_X$. Then we have
$$\phi_X(\sigma_X)=\frac{{\rm Sign}X-\deg X}{\chi(X)+\deg X-2(2-2g)}.$$
Here, ${\rm Sign}X$ is the signature of $X$ as a closed oriented 4-manifold and $\chi(X)$ is the Euler-Poincar\'e
characteristic of $X$, and $\deg X$ is the usual degree of $X$
(i.e., the number of intersecting points with a generic complementary dimensional plane to $X$).
\end{prop}

\begin{proof}
Let $L$ be a generic line of $\mathbb{P}_N^{\vee}$ as in the beginning of this subsection. As in \cite{La} (1.6.1)
the axis $A=\cap_{H\in L} H$ of the pencil meets $X$ transversally, and $i_L^*\mathcal{W}$ is the blow up of $X$
at the $\deg X$ points $A\cap X$ hence diffeomorphic to the connected sum
$X\# (\deg X)\overline{\mathbb{P}}_2$. Therefore we have
\begin{equation}
\label{eq:3-5-2}
{\rm Sign}(i_L^*\mathcal{W})={\rm Sign}X-\deg X
\end{equation}
and $\chi(i_L^*\mathcal{W})=\chi(X)+\deg X$.

Let $D_L=\{ H_1,\ldots,H_d\}$, where $d=\deg D_X$ and for $1\le i\le d$, let $\sigma_i\in \pi_1(U_L)$ be
the element obtained by substituting $H_i$ for $H^{\prime}$ in the definition of $\sigma$, see the beginning of
this subsection. As elements of $\pi_1(U^X)$, all $\sigma_i$ are lassos around $D_X$.

Let $D_i\subset L$ be a small closed 2-disk around $H_i$. We write by $f_L$ the projection
$i_L^*\mathcal{W}\rightarrow L$ and write $X_i=f_L^{-1}(D_i)$.
Let $X_0=X\setminus \coprod_i{\rm Int X_i}$. By Lemma \ref{lem:3-5-1}, $f_L^{-1}(H_i)$ is
irreducible hence the signature of $X_i$ is zero. Using the Novikov additivity, we have
${\rm Sign}(i_L^*\mathcal{W})={\rm Sign X_0}$.
By Meyer's signature formula (\cite{Meyer} Satz 1) and $\rho_X^*\tau_g=\delta \phi_X$, we have
\begin{equation}
\label{eq:3-5-3}
{\rm Sign}(i_L^*\mathcal{W})={\rm Sign X_0}=\sum_{i=1}^d \phi_X(\sigma_i)=d \phi_X(\sigma_X).
\end{equation}
On the other hand since all the singular fibers have one nodal singularity,
there are $d$ singular fibers with Euler contribution $+1$ (see \cite{BHPV}, (11.4) Proposition), thus
\begin{equation}
\label{eq:3-5-4}
d=\chi(i_L^*\mathcal{W})-2(2-2g)=\chi(X)+\deg X-2(2-2g).
\end{equation}
The proposition follows from (\ref{eq:3-5-2}), (\ref{eq:3-5-3}), and (\ref{eq:3-5-4}).
\end{proof}

Note that by (\ref{eq:3-5-4}) we can express $\deg D_X$ in terms of $\chi(X)$, $\deg X$, and $g$. The genus $g$ is
expressed as follows. For $H\in U^X$, let $C=H\cap X$. By the adjunction formula we have $c_1(C)=c_1(X)|_C-h$,
where $h$ is the hyperplane class, thus
$$2-2g=\chi(C)=\langle c_1(C),[C]\rangle=\langle c_1(X)h-h^2,[X] \rangle=\langle c_1(X)h,[X]\rangle- \deg X.$$

\begin{exple}{\rm
Let $m\ge 1$ and $n_1,\ldots,n_m \ge 2$ be integers and let $X\subset \mathbb{P}_{m+2}$
be a smooth complete intersection of type $(n_1,\ldots,n_m)$. Namely $X$ is given as the zero set of
some homogeneous polynomials $f_1,\ldots,f_m$ where $f_i$ is of degree $n_i$.
}
\end{exple}

\begin{prop}
\label{prop:3-5-3}
Let $X$ be as above and assume that $(m,n_1,\ldots,n_m)\neq (1,2)$. Then $D_X$
is a hypersurface. Let $\sigma_X\in \pi_1(U^X)$ be a lasso around $D_X$. We have
$$\phi_X(\sigma_X)=\frac{m-\displaystyle{\sum_{i=1}^mn_i^2}}
{3\left(\displaystyle{\frac{m^2+m}{2}}+\sum_{i=1}^mn_i^2-(m+1)\sum_{i=1}^mn_i+\sum_{i<j}n_in_j\right)}.$$
\end{prop}

\begin{proof}
We have $\deg X=n_1\cdots n_m$ and using the adjunction formula we can compute
\begin{equation}
\label{eq:3-5-5}
\chi(X)=c_2(X)=n_1\cdots n_m \left( \left( \begin{array}{c} m+3 \\ 2 \\ \end{array} \right)
+\sum_{i=1}^mn_i^2-(m+3)\sum_{i=1}^mn_i+\sum_{i<j}n_in_j \right),
\end{equation}
\begin{equation}
\label{eq:3-5-6}
{\rm Sign}X=\frac{n_1\cdots n_m}{3}\left(m+3-\sum_{i=1}^mn_i^2\right).
\end{equation}
For $H\in U^X$, $C=H\cap X$ is a smooth complete intersection of type $(n_1,\ldots,n_m,1)$.
Using the adjunction formula we have
$$2-2g=\chi(C)=n_1\cdots n_m\left(m+2-\sum_{i=1}^mn_i\right),$$
and following the argument in the proof of Proposition \ref{prop:3-5-2}, $\deg D_X$ is given by
$$\deg D_X=n_1\cdots n_m \left( \frac{m^2+m}{2}+\sum_{i=1}^mn_i^2-(m+1)\sum_{i=1}^mn_i+\sum_{i<j}n_in_j \right).$$
We claim that $\deg D_X$ is positive. If $m=1$, $\deg D_X=n_1(n_1-1)^2 >0$. If $m\ge 2$, Using the inequality
\begin{equation}
\label{eq:3-5-7}
\sum_{i=1}^mn_i^2 \ge \frac{2}{m-1}\sum_{i<j}n_in_j
\end{equation}
for $n_i\ge 0$ (this is easily derived from the geometric-arithmetic mean inequality), we have
\begin{eqnarray*}
\deg D_X &\ge& n_1\cdots n_m \left(\frac{m^2+m}{2}+\frac{m+1}{m-1}\sum_{i<j}n_in_j-(m+1)\sum_{i=1}^mn_i \right) \\
&=& n_1\cdots n_m \frac{m+1}{m-1}\sum_{i<j}(n_i-1)(n_j-1) \\
\end{eqnarray*}
Thus in any case $\deg D_X>0$, i.e., $D_X$ is a hypersurface. Also we can show $\chi(C)\le 0$ hence $g>0$ except
for the case $m=1$ and $n_1=2$. Finally, $X$ is simply connected. This follows from the Zariski theorem of
Lefschetz type, see \cite{La}, (8.1.1). Now Proposition \ref{prop:3-5-2} can be applied, and combining the above
computations all together we have the result.
\end{proof}

The next example is a generalization of the above, but it will illustrate that for a fixed variety,
how the value of the Meyer function on a lasso depends on a choice of its projective embedding.

\begin{exple}{\rm
Let $m\ge 0$, $n_1,\ldots,n_m \ge 2$, $n\ge 2$, and $d\ge 1$ be integers.
When $m=0$, we assume that $d$ $\ge 2$.
Let $v_d\colon \mathbb{P}_{m+n}\hookrightarrow \mathbb{P}_N$ be the Veronese embedding of degree $d$. Here,
$$N=\left( \begin{array}{c} n+m+d \\ d \\ \end{array}\right)-1.$$
Let $X$ be the $v_d$-image of a smooth complete intersection in $\mathbb{P}_{m+n}$ of type $(n_1,\ldots,n_m)$.
When $m=0$, $X$ is by definition the $v_d$-image of $\mathbb{P}_n$.
}
\end{exple}

\begin{prop}
\label{prop:3-5-4}
Let $X$ be as above and assume that $(d,m,n_1,\ldots,n_m)\neq (1,1,2)$ and $(n,d,m)\neq (2,2,0)$.
Then the $k$-th associated
variety $D_X$ is a hypersurface. Let $\sigma_X\in \pi_1(U^X)$ be a lasso around $D_X$. We have
$$\phi_X(\sigma_X)=\frac{\alpha_X}{\beta_X},$$
where
$$\alpha_X=\frac{m+n+1-\sum_{i=1}^mn_i^2-(n+1)d^2}{3},$$ and
\begin{eqnarray*}
\beta_X &=& \left( \begin{array}{c} m+n+1 \\ 2 \\ \end{array} \right)
+\displaystyle{\sum_{i=1}^m}n_i^2+\sum_{i<j}n_in_j-(m+n+1)\left(\sum_{i=1}^mn_i+nd\right) \\
& & +nd\sum_{i=1}^mn_i+\frac{(n^2+n)d^2}{2}.
\end{eqnarray*}
\end{prop}

\begin{proof}
Let $(N,M)\in \mathcal{L}_k(\mathbb{P}_N)$ be as in Lemma \ref{lem:3-4-3} and $X^{\prime}=M\cap X$. We may focus
on $X^{\prime}\subset M$. We will show that $D_{X^{\prime}}$ is a hypersurface and compute the value
$\phi_{X^{\prime}}(\sigma_{X^{\prime}})$, which must coincide with $\phi_X(\sigma_X)$, where
$\sigma_{X^{\prime}}\in \pi_1(U^{X^{\prime}})$ is a lasso around $D_{X^{\prime}}$.

First of all, the pull back $v_d^{-1}(X^{\prime})$ is a smooth complete intersection in $\mathbb{P}_{m+n}$
of type $(n_1,\ldots,n_m,\underbrace{d,\ldots,d}_{n-2})$. Thus $X^{\prime}$ is simply connected, and the invariants
$\chi(X^{\prime})$ and ${\rm Sign}X^{\prime}$ can be computed from (\ref{eq:3-5-5}), (\ref{eq:3-5-6}).
Also, we have $\deg X^{\prime}=\deg X=n_1\ldots n_md^n$. From these we can see that
${\rm Sign}X^{\prime}-\deg X^{\prime}={n_1\ldots n_md^{n-2}}\alpha_X$. For $W\in U^{X^{\prime}}\subset U^X$,
$C:=W\cap X=W\cap X^{\prime}$ is a smooth complete intersection in $\mathbb{P}_{m+n}$ of type
$(n_1,\ldots,n_m,\underbrace{d,\ldots,d}_{n-1})$. Thus the genus $g$ of $C$ is seen by 

$$2-2g=\chi(C)=n_1\ldots n_md^{n-1}\left(m+n+1-\sum_{i=1}^mn_i-(n-1)d\right).$$
It is easy to see that under the assumption, we have $\chi(C)\le 0$ hence $g>0$.
Using (\ref{eq:3-5-4}) and our knowledge of $\chi(X^{\prime})$, $\deg X^{\prime}$, and $\chi(C)$ gives
$$\deg D_{X^{\prime}}=n_1\ldots n_md^{n-2}\beta_X.$$
We claim that $\deg D_{X^{\prime}}$ is positive. Now we have the inequality
$$nd^2+\sum_{i=1}^mn_i^2\ge \frac{2}{m+n-1}\left(\sum_{i<j}n_in_j+nd\sum_{i=1}^mn_i
+\left( \begin{array}{c} n \\ 2 \\ \end{array} \right) d^2 \right),$$
the same kind of (\ref{eq:3-5-7}). Using this, we have
\begin{eqnarray*}
\deg D_{X^{\prime}} &\ge & n_1\ldots n_md^{n-2}\frac{m+n+1}{m+n-1}\left( \sum_{i<j}(n_i-1)(n_j-1)
+n\sum_{i=1}^m(n_i-1)(d-1) \right. \\
& &\left. +\left( \begin{array}{c} n \\ 2 \\ \end{array}\right)(d-1)^2 \right). \\
\end{eqnarray*}
This shows $\deg D_{X^{\prime}}>0$ except for the case $d=1$ and $m=1$. In this case, we have
$\deg D_{X^{\prime}}=n_1(n_1-1)^2>0$. Thus $D_{X^{\prime}}$ is a hypersurface, so is $D_X$.
Applying Proposition \ref{prop:3-5-2}, we have $\phi_X(\sigma_X)=\phi_{X^{\prime}}(\sigma_{X^{\prime}})
=({\rm Sign}X^{\prime}-\deg X^{\prime})/\deg D_{X^{\prime}}=\alpha_X/\beta_X$, as desired.
\end{proof}

\subsection{Bounded cohomology of $\pi_1(U^X)$}
For a group $G$, we denote by $H^*_b(G;\mathbb{R})$ the bounded cohomology group of $G$.
Namely, $H^*_b(G;\mathbb{R})$ is
the cohomology of the cochain complex of $\mathbb{R}$-valued {\it bounded} cochains of $G$.
In this subsection we show that the second bounded
cohomology of $\pi_1(U^X)$ is non-trivial under a certain mild condition.

\begin{prop}
\label{prop:3-6-1}
Let $X\subset \mathbb{P}_N$ be a smooth projective variety of dimension $\ge 2$
such that $D_X$ is a hypersurface. Suppose the value of $\phi_X$ on a lasso around $D_X$
is neither equal to $0$ nor $-1$. Then the bounded cohomology $H^2_b(\pi_1(U^X);\mathbb{R})$ is non-trivial and
the natural comparison map $H^2_b(\pi_1(U^X);\mathbb{R})\rightarrow H^2(\pi_1(U^X);\mathbb{R})$
is not injective.
\end{prop}

We need a lemma.

\begin{lem}
\label{lem:3-6-2}
Let $T\in \Gamma_g$ be the right hand Dehn twist along a non-separating simple closed curve on $\Sigma_g$.
Then for any integer $n\ge 1$, we have
$$\tau_g(T^{-1},T^{-n})=-1.$$
\end{lem}

\begin{proof}
We use the description (\ref{eq:1-3}) of $\tau_g$.
By the formulas (12) and (13) of \cite{Meyer}, it suffices to prove
${\rm Sign}(V_{A,A^n},\langle \ ,\ \rangle_{A,A^n})=-1$ where
$$A=\left( \begin{array}{cc} 1 & -1 \\ 0 & 1 \end{array} \right)$$
($A$ corresponds to the inverse of the right Dehn twist along
a non-separating simple closed curve on the torus). We have
$$V_{A,A^n}=\left\{ (x,y)\in \mathbb{R}^2\oplus \mathbb{R}^2; 
\left(\begin{array}{cc} 0 & 1 \\ 0 & 0 \end{array}\right)x+
\left(\begin{array}{cc} 0 & -n \\ 0 & 0 \end{array}\right)y=0 \right\},$$
thus the vectors $((1,0)^t,(0,0)^t)$, $((0,0)^t,(1,0)^t)$, and $((0,n)^t,(0,1)^t)$
form a basis for $V_{A,A^n}$. The presentation matrix of $\langle \ ,\ \rangle_{A,A^n}$
with respect to this basis is
$$\left( \begin{array}{ccc} 0 & 0 & 0 \\ 0 & 0 & 0 \\ 0 & 0 & -n(n+1) \\ \end{array} \right).$$
This completes the proof.
\end{proof}

\begin{proof}[Proof of Proposition \ref{prop:3-6-1}]
Note that $\tau_g$ is a bounded 2-cocycle of $\Gamma_g$. More precisely,
for $f_1,f_2 \in \Gamma_g$ we have $|\tau_g(f_1,f_2)|\le 4g$ .
Thus $\rho_X^*\tau_g$ is also a bounded 2-cocycle. Since $\phi_X$ is a unique 1-cochain
cobounding $\rho_X^*\tau_g$, it suffices to show that $\phi_X$ is unbounded.

Let $\sigma_X$ be a lasso around $D_X$. By the Picard-Lefschetz formula,
$\rho_X(\sigma_X)$ is the inverse of the right hand Dehn twist along a simple closed curve.
We claim that this curve is non-separating.
For, if this is separating, $\rho_X(\sigma_X)\in \Gamma_g$ does act trivially
on the homology of $\Sigma_g$. Combining this with the fact that
$\pi_1(U^X)$ is normally generated by $\sigma_X$, we deduce that
the image $\rho_X(\pi_1(U^X))$ acts trivially on the homology of $\Sigma_g$.
Hence $\rho_X^*\tau_g$ is zero as a cocycle.
Since $\delta \phi_X=\rho_X^*\tau_g=0$ and $\phi_X(\sigma_X)\neq 0$ it follows
that $\phi_X$ is a non-trivial homomorphism from $\pi_1(U^X)$ to $\mathbb{Q}$,
contradicting to Lemma \ref{lem:3-3-1}.

Now by $\delta\phi_X=\rho_X^*\tau_g$ and Lemma \ref{lem:3-6-2}, we have
$$\phi_X(\sigma_X^n)=n\phi_X(\sigma_X)-\sum_{i=1}^{n-1}\tau_g(\rho_X(\sigma_X),\rho_X(\sigma_X^i))
=n\phi_X(\sigma_X)+n-1$$
for $n\ge 1$. Since $\phi_X(\sigma_X)\neq -1$, this shows the unboundedness of $\phi_X$.
\end{proof}

It is known that if a discrete group is amenable, then its bounded cohomology vanishes in positive degrees
(see \cite{Grom}). Thus:

\begin{cor}
\label{cor:3-6-3}
Let $X\subset \mathbb{P}_N$ be a smooth projective variety satisfying the hypothesis of
Proposition \ref{prop:3-6-1}. Then the fundamental group $\pi_1(U^X)$ is not amenable.
\end{cor}

As an example, when $X$ is the one of those in Proposition \ref{prop:3-5-4},
we can check that $\alpha_X < 0$ and $\alpha_X+\beta_X > 0$. Therefore
Proposition \ref{prop:3-6-1} can be applied to this situation.

\section{Applications to local signatures}
\subsection{An approach to local signatures via Meyer functions}
Let $\mathbb{M}_g$ be the moduli space of compact Riemann surfaces of genus $g$ and $\mathcal{A}$ a subset of 
$\mathbb{M}_g$. We introduce the notion of an $\mathcal{A}$-fibration and a local signature
with respect to $\mathcal{A}$.

\begin{dfn}
\label{dfn:4-1-1}
Let $B$ be a topological space.
\begin{enumerate}
\item
A triple $\xi=(\mathcal{C},p,B)$ is called an \textit{$\mathcal{A}$-family} (on $B$)
if $p\colon \mathcal{C}\rightarrow B$
is a continuous family of compact Riemann surfaces with each fiber being an element of $\mathcal{A}$.
\item
Let $\xi_0,\xi_1$ be $\mathcal{A}$-families on $B$. They are called \textit{isotopic} if there exists an
$\mathcal{A}$-family $\xi$ on $B\times [0,1]$ such that for $i=0,1$, the restriction of $\xi$ to
$B\times \{ i\}$ is isomorphic to $\xi_i$ as continuous family of Riemann surfaces on $B$,
where $B\times \{i\}$ is identified with $B$ by $(b,i)\mapsto b$.
\item
We denote by $\mathcal{A}(B)$ the set of isotopy classes of $\mathcal{A}$-families over $B$.
\end{enumerate}
\end{dfn}

Let $\xi=(\mathcal{C},p,B)$ be an $\mathcal{A}$-family and $\psi\colon B^{\prime}\rightarrow B$ a continuous map.
By taking the fiber product of $\psi$ and $p$, the pull back $\psi^*\xi$ by $\psi$ is naturally defined as
an $\mathcal{A}$-family on $B^{\prime}$. In this way we get a category of $\mathcal{A}$-families, which
we denote by $\mathcal{A}^{\prime}$. Moreover this association induces the map
$\eta_{\xi,B^{\prime}} \colon [B^{\prime},B]\rightarrow \mathcal{A}(B^{\prime})$, where $[B^{\prime},B]$
is the set of homotopy classes of continuous maps from $B^{\prime}$ to $B$.

\begin{dfn}
\label{dfn:4-1-2}
An $\mathcal{A}$-family $\xi_u=(\mathcal{C}_u,p_u,B_u)$ over a path connected space $B_u$ having the homotopy type of a 
CW-complex is called \textit{universal} if $\eta_{\xi_u,B}$ is bijective for any topological space $B$ having the
homotopy type of a CW-complex. We denote by $\rho_u\colon \pi_1(B_u)\rightarrow \Gamma_g$ the topological monodromy of
$p_u\colon \mathcal{C}_u\rightarrow B_u$.
\end{dfn}

A universal $\mathcal{A}$-family is uniquely determined up to isotopy if it exists:
if $\xi_u^{\prime}=(\mathcal{C}_u^{\prime},p_u^{\prime},B_u^{\prime})$ is another universal $\mathcal{A}$-family
then there exist continuous maps $\psi\colon B_u \rightarrow B_u^{\prime}$ and
$\psi^{\prime}\colon B_u^{\prime} \rightarrow B_u$ such that
$\psi^*\xi_u^{\prime}$ and $\xi_u$ (resp. ${\psi^{\prime}}^*\xi_u$ and $\xi_u^{\prime}$) are isotopic.
In particular $B_u$ and $B_u^{\prime}$ are homotopy equivalent.

In some situations as we will see, we can construct a universal $\mathcal{A}$-family from a certain
$\mathcal{A}$-family with group action. The following proposition is used to verify the universality
of such a family. To state the proposition we prepare a terminology.
Let $\xi^0=(\mathcal{C}^0,p^0,B^0)$ be a $\mathcal{A}$-family on a connected $C^{\infty}$-manifold $B^0$, and
let $\mathcal{G}$ be a Lie group acting on $\mathcal{C}^0$ and $B^0$ from the left such
that $p^0$ is $\mathcal{G}$-equivariant. Let $\mathcal{G}^{\prime}$ be the category defined as follows:
the objects consist of $(P,\pi,B,E)$ such that $\pi\colon P\rightarrow B$ is a principal
$\mathcal{G}$-bundle on $B$ (the $\mathcal{G}$-action on $P$ being from the left)
and $E\colon P\rightarrow B^0$ is a $\mathcal{G}$-equivariant map, and the morphisms from
$(P,\pi,B,E)$ and $(P^{\prime},\pi^{\prime},B^{\prime},E^{\prime})$ are the bundle maps
from $P$ to $P^{\prime}$ compatible with $E$ and $E^{\prime}$.

\begin{prop}[A criterion for universality]
\label{prop:4-1-3}
Let $\mathcal{A}$ be a subset of $\mathbb{M}_g$ and $\xi^0=(\mathcal{C}^0,p^0,B^0)$ an $\mathcal{A}$-family
as above. Suppose there is a covariant functor from $\mathcal{A}^{\prime}$ to $\mathcal{G}^{\prime}$
associating an $\mathcal{A}$-family $\xi=(\mathcal{C},p,B)$ with $(P(\xi),\pi,B,E_{\xi})$, and
satisfying the following conditions:
\begin{enumerate}
\item
${E_{\xi}}^*\xi^0$ and $\pi^*\xi$ are isomorphic as continuous families of Riemann surfaces on $P(\xi)$.
\item
As for the object associated to $\xi^0$, we can take a trivial $\mathcal{G}$-bundle 
$P(\xi^0)=\mathcal{G}\times B^0$ and a $\mathcal{G}$-equivariant map $E_{\xi^0}$ such that $E_{\xi^0}(g,b)=g\cdot b$.
Moreover, for any $g\in \mathcal{G}$ the bundle map $\bar{g}\colon P(\xi^0)\rightarrow P(\xi^0)$
induced by the maps $B^0\rightarrow B^0$, $b\mapsto g\cdot b$ and
$\mathcal{C}^0\rightarrow \mathcal{C}^0$, $c\mapsto g\cdot c$ is given by
$\bar{g}(g^{\prime},b)=(g^{\prime}g^{-1},g\cdot b)$.
\end{enumerate}
Let $E\mathcal{G}\rightarrow B\mathcal{G}$ be a universal principal $\mathcal{G}$-bundle
(the $\mathcal{G}$-action on $E\mathcal{G}$ being from the right).
Taking the Borel constructions $B^0_{\mathcal{G}}=E\mathcal{G}\times_{\mathcal{G}}B^0$ and 
$\mathcal{C}^0_{\mathcal{G}}$ we obtain an $\mathcal{A}$-family 
$\xi^0_{\mathcal{G}}=(\mathcal{C}^0_{\mathcal{G}},p^0_{\mathcal{G}},B^0_{\mathcal{G}})$.
Then, $\xi^0_{\mathcal{G}}$ is a universal $\mathcal{A}$-family.
\end{prop}
\begin{proof}
Let $B$ be a space having the homotopy type of a CW-complex. For simplicity we write 
$\eta=\eta_{\xi^0_{\mathcal{G}},B}$. We construct a candidate for the inverse of $\eta$.
Let $\xi=(\mathcal{C},p,B)$ be an $\mathcal{A}$-family on $B$. Take a principal $\mathcal{G}$-bundle $P=P(\xi)$
and a $\mathcal{G}$-equivariant map $E=E_{\xi}$ associated to $\xi$. Considering the Borel construction
$P_{\mathcal{G}}=E\mathcal{G}\times_{\mathcal{G}}P$, let $T\colon P_{\mathcal{G}}\rightarrow B$
be the map induced from the projection $\pi\colon P\rightarrow B$. This is an $E\mathcal{G}$-bundle, thus by Dold's 
theorem it has a section: a map $\zeta\colon B\rightarrow P_{\mathcal{G}}$ such that $T\circ \zeta={\rm id}_B$.
Let $E_{\mathcal{G}}\colon P_{\mathcal{G}}\rightarrow B^0_{\mathcal{G}}$ be the map induced from $E$.
Now the isomorphism $E^*\xi^0\cong \pi^*\xi$ induces the isomorphism
${E_{\mathcal{G}}}^*\xi^0_{\mathcal{G}}\cong T^*\xi$ and $E_{\mathcal{G}}\circ \zeta$ is a continuous map from
$B$ to $B^0_{\mathcal{G}}$ such that
$$(E_{\mathcal{G}}\circ \zeta)^*\xi^0_{\mathcal{G}}=\zeta^*{E_{\mathcal{G}}}^*\xi^0_{\mathcal{G}}\cong
\zeta^*T^*\xi=(T\circ \zeta)^*\xi=\xi.$$
This shows $\eta$ is surjective. In fact, using the functoriality we can show that the homotopy class of
$E_{\mathcal{G}}\circ \zeta$ depends only on the isotopy class of $\xi$. In this way we have the map
$\theta\colon \mathcal{A}(B)\rightarrow [B,B^0_{\mathcal{G}}]$ satisfying
$\eta \circ \theta={\rm id}_{\mathcal{A}(B)}$.

Here we consider the above construction applied to $\xi^0_{\mathcal{G}}$.
We have $P(\xi^0)=\mathcal{G}\times B^0$ with the projection $\pi^0\colon P(\xi^0)\rightarrow B^0$,
$(g,b)\mapsto b$ and the $\mathcal{G}$-equivariant map $E_{\xi^0}\colon P(\xi^0)\rightarrow B^0$,
$(g,b)\mapsto g\cdot b$. By the functoriality, the $\mathcal{G}$-bundle
$\pi_u\colon  P^u\rightarrow B^0_{\mathcal{G}}$ and the $\mathcal{G}$-equivariant map
$E^u\colon P^u\rightarrow B^0$ associated to $\xi^0_{\mathcal{G}}$ is described as follows.

Take the Borel construction $P^u:=E\mathcal{G}\times_{\mathcal{G}}P(\xi^0)$ where the $\mathcal{G}$-action on
$P(\xi^0)$ is given by $g\cdot (g^{\prime},b)=(g^{\prime}g^{-1},g\cdot b)$. Define
$\pi_u \colon P^u\rightarrow B^0_{\mathcal{G}}$ and $E^u\colon P^u\rightarrow B^0$ by
$$\pi_u([e,(g,b)])=[e,b], \textrm{ and } E^u([e,(g,b)])=g\cdot b,$$
where $[e,b]$ denotes the element of $B^0_{\mathcal{G}}$ represented by $(e,b)\in E\mathcal{G}\times B^0$, etc.
Given the $\mathcal{G}$-action on $P^u$ by $g\cdot [e,(g^{\prime},b)]=[e,(gg^{\prime},b)]$,
$\pi_u$ is a principal $\mathcal{G}$-bundle and $E^u$ is $\mathcal{G}$-equivariant. Also the isomorphism
${E_{\xi^0}}^*\xi_0\cong {\pi^0}^*\xi_0$ induces ${E^u}^*\xi_0 \cong {\pi_u}^*\xi^0_{\mathcal{G}}$.
Notice that $T_u\colon P^u_{\mathcal{G}}\rightarrow B^0_{\mathcal{G}}$ has a section $\zeta_u$ given by
$\zeta_u([e,b])=[e,[e,({\rm id}_{\mathcal{G}},b)]]$.

Now we show $\theta \circ \eta={\rm id}_{[B,B^0_{\mathcal{G}}]}$, which will complete the proof.
Let $\psi\colon B\rightarrow B^0_{\mathcal{G}}$ be a continuous map.
By the functoriality, we can use the fiber product $\psi^*P^u$ as the $\mathcal{G}$-bundle associated to
$\psi^*\xi^0_{\mathcal{G}}$. Pulling back $\zeta_u$, we have a section
$\psi^*\zeta_u$ of $T\colon (\psi^*P^u)_{\mathcal{G}}\rightarrow B$ which makes the following diagram commutative.
$$\xymatrix{ (\psi^*P^u)_{\mathcal{G}} \ar[r]^{\bar{\psi}_{\mathcal{G}}} & P^u_{\mathcal{G}} \ar[r]^{E^u_{\mathcal{G}}}
& B^0_{\mathcal{G}} \\ B \ar[r]_{\psi} \ar[u]^{\psi^*\zeta_u} & B^0_{\mathcal{G}} \ar[u]^{\zeta_u} & \\ }$$
Notice that $E^u_{\mathcal{G}} \circ \zeta_u={\rm id}_{B^0_{\mathcal{G}}}$.
Following the construction of $\theta$ we have
$$\theta([\psi^*\xi^0_{\mathcal{G}}])=[E^u_{\mathcal{G}} \circ \bar{\psi}_{\mathcal{G}} \circ \psi^*\zeta_u]
=[E^u_{\mathcal{G}} \circ \zeta_u \circ \psi]=[\psi],$$
which shows $\theta \circ \eta={\rm id}_{[B,B^0_{\mathcal{G}}]}$.
\end{proof}

\begin{dfn}
\label{dfn:4-1-4}
\begin{enumerate}
\item
Let $\Delta$ be a closed oriented 2-disk with the center $b$. A 4-tuple $\mathcal{F}=(S,f,\Delta,b)$ is called an
\textit{$\mathcal{A}$-degeneration} if $S$ is a $C^{\infty}$-manifold of dimension 4 and
$f \colon S\rightarrow \Delta$
is a proper surjectice $C^{\infty}$-map, and the restriction of $f$ to $\Delta \setminus \{ b\}$ is given
a structure of $\mathcal{A}$-family. We denote by $\xi_{\mathcal{F}}$ this $\mathcal{A}$-family.
\item
Let $\mathcal{F}=(S,f,\Delta,b)$ and $\mathcal{F}^{\prime}=(S^{\prime},f^{\prime},\Delta^{\prime},b^{\prime})$ 
be $\mathcal{A}$-degenerations. We say $\mathcal{F}$ and $\mathcal{F}^{\prime}$ are \textit{equivalent} if taking
suitably smaller disks $\Delta_0\subset \Delta$ with $b\in \Delta_0$ and $\Delta_0^{\prime}\subset \Delta^{\prime}$
with $b^{\prime} \in \Delta_0^{\prime}$, there exist an orientation preserving homeomorphism
$\psi\colon (\Delta_0,b)\rightarrow (\Delta_0^{\prime},b^{\prime})$ such that
$\psi^*\xi_{\mathcal{F}^{\prime}}$ is isotopic to the restriction of $\xi_{\mathcal{F}}$
to $\Delta_0 \setminus \{b\}$.

\item
We denote by $\mathcal{A}_{loc}$ the set of all equivalence classes of $\mathcal{A}$-degenerations.
We often identify an element of $\mathcal{A}_{loc}$ with its representative.
Each element of $\mathcal{A}_{loc}$ is called a fiber germ. A \textit{smooth fiber germ}
is an element of $\mathcal{A}_{loc}$ obtained by an $\mathcal{A}$-family $\xi=(\mathcal{C},p,\Delta)$.
\end{enumerate}
\end{dfn}

\begin{dfn}
\label{dfn:4-1-5}
Let $M$ be a closed oriented 4-manifold and $B$ a closed oriented 2-manifold.
A proper surjective $C^{\infty}$-map $f\colon M\rightarrow B$ is called an \textit{$\mathcal{A}$-fibration} if
there exist finitely many points $b_1,\ldots,b_m \in B$ such that the restriction of $f$ to
$B\setminus \{ b_1,\ldots,b_m\}$ is given a structure of $\mathcal{A}$-faimly.
\end{dfn}

The triple $(M,f,B)$ is a fibered 4-manifold in the sense of section 1.
For an $\mathcal{A}$-fibration $f\colon M\rightarrow B$, let $\mathcal{F}_i$ be the element of $\mathcal{A}_{loc}$
obtained by restricting $f$ to a small closed disk neighborhood $\Delta_i$ of $b_i$. We formulate
the notion of a local signature in our setting.

\begin{dfn}
\label{dfn:4-1-6}
Let $\mathcal{A}$ be a subset of $\mathbb{M}_g$. A function
$\sigma_{\mathcal{A}}\colon \mathcal{A}_{loc}\rightarrow \mathbb{Q}$ is called a local signature
with respect to $\mathcal{A}$ if
\begin{enumerate}
\item
for a smooth fiber germ $\mathcal{F}$, $\sigma_{\mathcal{A}}(\mathcal{F})=0$, and
\item
for any $\mathcal{A}$-fibration $f\colon M\rightarrow B$, we have \textit{the global signature formula}:
\begin{equation}
\label{eq:4-1-1}
{\rm Sign}(M)=\sum_{i=1}^m \sigma_{\mathcal{A}}(\mathcal{F}_i).
\end{equation}
\end{enumerate}
\end{dfn}

\begin{prop}
\label{prop:4-1-7}
Let $\mathcal{A}$ be a subset of $\mathbb{M}_g$. Suppose there exist a universal $\mathcal{A}$-family
$\xi_u=(\mathcal{C}_u,p_u,B_u)$ and a $\mathbb{Q}$-valued 1-cochain
$\phi_{\mathcal{A}}\colon \pi_1(B_u)\rightarrow \mathbb{Q}$ 
such that $\delta \phi_{\mathcal{A}}=\rho_u^*\tau_g$. Then there exists a local signature with respect to $\mathcal{A}$.
\end{prop}

\begin{proof}
Let $\mathcal{F}=(S,f,\Delta,b)\in \mathcal{A}_{loc}$.
Since $\eta=\eta_{\xi_u,\Delta \setminus \{ b\}}$ is bijective, there exists
uniquely up to homotopy a continuous map $g_{\mathcal{F}}\colon \Delta \setminus \{ b\} \rightarrow B_u$ such that
$\eta([g_{\mathcal{F}}])=[\xi_{\mathcal{F}}]$.
We denote by $\partial \Delta$ the element of $\pi_1(\Delta \setminus \{ b\})$
represented by the loop going once around the boundary of $\Delta$ by counter clockwise manner.
Then we obtain an element
$x_{\mathcal{F}}={g_{\mathcal{F}}}_*(\partial \Delta)\in \pi_1(B_u)$, which is uniquely determined up to conjugacy.
Since the equality $\delta \phi_{\mathcal{A}}=\rho_u^*\tau_g$ implies that $\phi_{\mathcal{A}}$ is a class function
(see section 1), the value $\phi_{\mathcal{A}}(x_{\mathcal{F}})$ is well defined.

Now define $\sigma_{\mathcal{A}}\colon \mathcal{A}_{loc}\rightarrow \mathbb{Q}$ by
\begin{equation}
\label{eq:4-1-2}
\sigma_{\mathcal{A}}(\mathcal{F})=\phi_{\mathcal{A}}(x_{\mathcal{F}})+{\rm Sign}(S).
\end{equation}
If $\mathcal{F}$ is a smooth fiber germ, $g_{\mathcal{F}}$ extends to a continuous map from $\Delta$.
So $x_{\mathcal{F}}\in \pi_1(B_u)$ is trivial, hence $\phi_{\mathcal{A}}(x_{\mathcal{F}})=0$.
Also we have ${\rm Sign}(S)=0$ since topologically $f\colon S\rightarrow \Delta$ is just a trivial $\Sigma_g$-bundle.
The first condition in Definition \ref{dfn:4-1-6} is verified. The second condition is verified by an argument
similar to the proof of Proposition \ref{prop:3-5-2}, so we omit the detail (see also \cite{Ku}, Theorem 7.2).
\end{proof}

\begin{dfn}
\label{dfn:4-1-8}
For $\mathcal{F}\in \mathcal{A}_{loc}$, we call $x_{\mathcal{F}}\in \pi_1(B_u)$
appeared in the proof of Proposition \ref{prop:4-1-7} the \textit{lifted monodromy}. This is
uniquely determined up to conjugacy.
\end{dfn}

\subsection{Fibrations of rank 4 non-hyperelliptic curves of genus 4}
Let $C$ be a non-hyperelliptic Riemann surface of genus 4. Its canonical image is a
$(2,3)$ complete intersection in $\mathbb{P}_3$ hence is contained
in a uniquely determined quadric. We say $C$ is \textit{of rank 4} if this quadric is of rank 4.
Let $\mathcal{R}^4\subset \mathbb{M}_4$ be the set of rank 4 non-hyperelliptic Riemann surfaces
of genus 4. $\mathcal{R}^4$ is Zariski open in $\mathbb{M}_4$.

Let $s \colon \mathbb{P}_1\times \mathbb{P}_1 \rightarrow \mathbb{P}_3$ be the Segre embedding.
Explicitly, $s$ is given by
$$s([a_0:a_1],[b_0:b_1])=[a_0b_0:a_0b_1:a_1b_0:a_1b_1],$$
using the homogeneous coordinates.
Let $V_{3,3}=\mathbb{C}[a_0,a_1]^3 \otimes \mathbb{C}[b_0,b_1]^3$
be the space of $(3,3)$ homogeneous polynomials, and let
$$s_{3,3}\colon \mathbb{P}_1\times \mathbb{P}_1 \rightarrow
\mathbb{P}({V_{3,3}}^{\vee})\cong \mathbb{P}(V_{3,3})^{\vee} \cong \mathbb{P}_{15}$$
be the embedding induced from the evaluation map
$\mathbb{C}^2 \times \mathbb{C}^2 \rightarrow {V_{3,3}}^{\vee}={\rm Hom}(V_{3,3},\mathbb{C})$.
Set
$$X={\rm Im}(s_{3,3}).$$

Consider the group $\mathcal{G}={\rm Aut}(\mathbb{P}_1\times \mathbb{P}_1)$.
Of course $\mathcal{G}$ acts on $\mathbb{P}_1\times \mathbb{P}_1$ (from the left),
inducing an action of
$\mathcal{G}$ on $\mathbb{P}_{15}$ so that $s_{3,3}$ is $\mathcal{G}$-equivariant.
Moreover $\mathcal{G}$
naturally acts on $\mathbb{P}(V_{3,3})$ from the left.

Let $D_X\subset \mathbb{P}_{15}^{\vee}=\mathbb{P}(V_{3,3})$ be the dual variety of $X$ and
$U^X=\mathbb{P}(V_{3,3})\setminus D_X$. $D_X$ is preserved by the $\mathcal{G}$-action.
Also $\mathcal{G}$ acts on $\mathcal{C}^X\subset \mathbb{P}_{15}\times U^X$ diagonally, and
the projection $p_X\colon \mathcal{C}^X\rightarrow U^X$ is $\mathcal{G}$-equivariant.
Note that for $W\in U^X$, the fiber $p_X^{-1}(W)$ is isomorphic to the smooth curve in
$\mathbb{P}_1\times \mathbb{P}_1$ determined by a $(3,3)$ homogeneous polynomial,
which is an element of $\mathcal{R}^4$ since the restriction of $s$ to the curve gives its canonical embedding
and the image is contained in $s(\mathbb{P}_1\times \mathbb{P}_1)$, which is a smooth quadric
$x_0x_3-x_1x_2=0$. Thus, $\xi^X=(\mathcal{C}^X,p_X,U^X)$ is a $\mathcal{R}^4$-family.

Now we will show that $\xi^X$ and the $\mathcal{G}$-action on it satisfies the conditions in Proposition
\ref{prop:4-1-3}. We need to consrtuct a principal $\mathcal{G}$-bundle from a $\mathcal{R}^4$-family.

First we consider the case of a single element $C\in \mathcal{R}^4$. We denote by
$\Omega^1(C)$ the space of holomorphic 1-forms on $C$.
The unique quadric containing the canonical image of $C$ corresponds to the 1-dimensional kernel of the natural map
$t_2\colon {\rm Sym}^2\Omega^1(C)\rightarrow H^0(C;K_C^{\otimes 2})$. Here, $K_C$ is the canonical bundle of $C$.
Note that $t_2$ is surjective by Max Noether's theorem (see \cite{GH}, p.\ 253).

If we take a basis $\omega=(\omega_0,\omega_1,\omega_2,\omega_3)$ of $\Omega^1(C)$, an explicit form of
${\rm Ker}(t_2)$ is obtained as follows. Let $\varphi_0,\varphi_1,\varphi_2,\varphi_3\in \Omega^1(C)^*$
be the dual basis of $\omega$. Then ${\rm Sym}^2\Omega^1(C)$ is identified with the space $S_4$ of
$4\times 4$ symmetric matrices by assigning $B\in S_4$ with the quadratic function
$$\Omega^1(C)^* \rightarrow \mathbb{C}, \ x_0\varphi_0+x_1\varphi_1+x_2\varphi_2+x_3\varphi_3
\mapsto (x_0,x_1,x_2,x_3)B(x_0,x_1,x_2,x_3)^t.$$
Hence a choice of a basis $\omega$ of $\Omega^1(C)$ determines the element $B(\omega)\in \mathbb{P}(S_4)$
corresponding to ${\rm Ker}(t_2)$, and the image of the canonical map
$\iota_{\omega}\colon C\hookrightarrow \mathbb{P}_3, \ c\mapsto [\omega_0(c):\omega_1(c):\omega_2(c):\omega_3(c)]$
is contained in the quadric determined by $B(\omega)$.

Let $P(C)$ be the set of $\omega$ modulo $\mathbb{C}^*$ such that
the quadric determined by $B(\omega)$ is equal to $\{ x_0x_3-x_1x_2=0\}$. Namely,
$$P(C)=\{ \omega \ {\rm mod}\mathbb{C}^*; \omega \textrm{ is a basis of } \Omega^1(C) \textrm{ and }
B(\omega)=H \},$$
where $H\in \mathbb{P}(S_4)$ is represented by
$$\left( \begin{array}{cccc}
0 & 0 & 0 & 1 \\
0 & 0 & -1 & 0 \\
0 & -1 & 0 & 0 \\
1 & 0 & 0 & 0 \\
\end{array} \right).$$
Now consider the group $PO_4^H(\mathbb{C})=\{ A\in PGL(4); A^tHA=H \}$, which acts on $P(C)$ from the left
freely and transitively by
$A\cdot (\omega_0,\omega_1,\omega_2,\omega_3) \ {\rm mod}\mathbb{C}^*=(\omega_0,\omega_1,\omega_2,\omega_3)
A^t \ {\rm mod}\mathbb{C}^*$. In fact, this group is isomorphic to $\mathcal{G}$ and the isomorphism is induced
by the action of $PO_4^H(\mathbb{C})$ on $\mathbb{P}_3$ (as a subgroup of $PGL(4)$) preserving
$s(\mathbb{P}_1\times \mathbb{P}_1)$. Therefore, $\mathcal{G}$ acts on $P(C)$ freely and transitively.
Finally, define the map
$$E_C\colon P(C)\rightarrow U^X$$
as follows. Again by Max Noether's theorem, the natural map
$t_3\colon {\rm Sym}^3\Omega^1(C)\rightarrow H^0(C;K_C^{\otimes 3})$ is surjective.
Choose $h\in {\rm Ker}(t_3)$ which is not divided by elements of ${\rm Ker}(t_2)$.
Let $\omega\in P(C)$. Then $h$ is identified with a homogeneous polynomial of degree 3 in
determinates $x_0,x_1,x_2,x_3$. We denote it by $h^{\omega}$.
The canonical image $\iota_{\omega}(C)$ is given by $x_0x_3-x_1x_2=h^{\omega}(x_0,x_1,x_2,x_3)=0$.
Set
$$E_C(\omega)=h^{\omega}(a_0b_0,a_0b_1,a_1b_0,a_1b_1)\in \mathbb{P}(V_{3,3}).$$
Since the zero set of $h^{\omega}(a_0b_0,a_0b_1,a_1b_0,a_1b_1)$ is isomorphic to $C$, we have $E_C(\omega)\in U^X$.
$E_C(\omega)$ does not depend on the choice of $h$, and we can verify that $E_C$ is $\mathcal{G}$-equivariant.

Now let $\xi=(\mathcal{C},p,B)$ be a $\mathcal{R}^4$-family. Applying the above construction to
all the fibers, we get a principal $\mathcal{G}$-bundle
$$P(\xi)=\bigcup_{b\in B}P(p^{-1}(b))$$
and by piecing together $E_{p^{-1}(b)}$, $b\in B$, we get a $\mathcal{G}$-equivariant map
$$E_{\xi}\colon P(\xi)\rightarrow U^X.$$
The first condition in Proposition \ref{prop:4-1-3} is clear from the construction. So far we have only
used the objects arising from holomorphic 1-forms on Riemann surfaces, which behave naturally
under pull back by biholomorphic maps. Thus the functoriality is also true. As to the second condition,
we can describe $P(\xi^X)$ as follows. For $W\in U^X$, $W\cap X$ is isomorphic to
$v_3^{-1}(W\cap X)\subset \mathbb{P}_3$.
Let $\omega_W$ be the basis of $\Omega^1(W\cap X)$ corresponding to the homogeneous coordinates
$[x_0:x_1:x_2:x_3]$ of $\mathbb{P}_3$. Then the isomorphism
$$P(\xi_X)\cong \mathcal{G}\times U^X$$
is given by assigning $(A,W)\in \mathcal{G}\times U^X$ with $A\cdot \omega_W$. We can check that for $A\in \mathcal{G}$
and $W\in U^X$, $\omega_W$ corresponds to $A^{-1}\cdot \omega_{A\cdot W}$ by the isomorphism
$W\cap X\rightarrow (A\cdot W)\cap X$ induced by $A$. The second condition follows from this.
Applying Proposition \ref{prop:4-1-3}, we have a universal $\mathcal{R}^4$-family
$\xi^X_{\mathcal{G}}=(\mathcal{C}^X_{\mathcal{G}},p_u,U^X_{\mathcal{G}})$. Here, $p_u=(p_X)_{\mathcal{G}}$.

\begin{thm}
\label{thm:4-2-1}
Let $\mathcal{R}^4$ be the set of rank 4 non-hyperelliptic Riemann surfaces of genus 4 and $X$, $\mathcal{G}$
as above. Then $\xi^X_{\mathcal{G}}=(\mathcal{C}^X_{\mathcal{G}},p_u,U^X_{\mathcal{G}})$ is a
universal $\mathcal{R}^4$-family. We denote by $\rho_u\colon \pi_1(U^X_{\mathcal{G}})\rightarrow \Gamma_4$
the topological monodromy of $p_u\colon \mathcal{C}^X_{\mathcal{G}}\rightarrow U^X_{\mathcal{G}}$.
Then there exists a unique $\mathbb{Q}$-valued 1-cochain
$\phi_{\mathcal{R}^4}\colon \pi_1(U^X_{\mathcal{G}}) \rightarrow \mathbb{Q}$ whose coboundary equals
to $\rho_u^*\tau_4$.
\end{thm}
\begin{proof}
We only have to prove the latter part. Consider the map
$U^X_{\mathcal{G}}\rightarrow B\mathcal{G}$ induced from the projection $E\mathcal{G}\rightarrow B\mathcal{G}$.
This is a $U^X$-bundle. By the homotopy exact sequence, we have the exact sequence
\begin{equation}
\label{eq:4-2-1}
\pi_1(\mathcal{G})\rightarrow \pi_1(U^X)\stackrel{i}{\rightarrow} \pi_1(U^X_{\mathcal{G}})
\rightarrow \pi_0(\mathcal{G})\rightarrow *,
\end{equation}
where $i$ is induced from the inclusion.
But $\pi_1(\mathcal{G})\cong \mathbb{Z}/2\mathbb{Z} \oplus \mathbb{Z}/2\mathbb{Z}$ and
$\pi_0(\mathcal{G})$ consists of two points. This follows from the fact $\mathcal{G}$ is isomorphic to
the semi direct product $(PGL(2)\times PGL(2))\ltimes \mathbb{Z}/2\mathbb{Z}$. So (\ref{eq:4-2-1})
shows that $i^*\colon H^2(\pi_1(U^X_{\mathcal{G}});\mathbb{Q})\rightarrow H^2(\pi_1(U^X);\mathbb{Q})$
is injective. By Theorem \ref{thm:3-1-1}, we have $i^*\rho_u^*[\tau_4]=0\in H^2(\pi_1(U^X);\mathbb{Q})$.
Therefore we also have $\rho_u^*[\tau_4]=0\in H^2(\pi_1(U^X_{\mathcal{G}});\mathbb{Q})$.
This shows the existence of $\phi_{\mathcal{R}^4}$.

On the other hand, in the proof of Lemma \ref{lem:3-3-1} we have seen that the abelianization of $\pi_1(U^X)$
is a cyclic group of finite order. Combining this fact with (\ref{eq:4-2-1}), we see that the abelianization
of $\pi_1(U^X_{\mathcal{G}})$ is a finite abelian group. This shows the uniqueness of $\phi_{\mathcal{R}^4}$.
\end{proof}

Combining this with Proposition \ref{prop:4-1-7}, we have

\begin{cor}
\label{cor:4-2-2}
Let $\mathcal{R}^4$ be the set of rank 4 non-hyperelliptic Riemann surfaces of genus 4.
Then the formula
\begin{equation}
\label{eq:4-2-2}
\sigma_{\mathcal{R}^4}(\mathcal{F})=\phi_{\mathcal{R}^4}(x_{\mathcal{F}})+{\rm Sign}(S)
\end{equation}
for $\mathcal{F}=(S,\pi,\Delta,b)\in \mathcal{R}^4_{loc}$ (see (\ref{eq:4-1-2})) gives
a local signature with respect to $\mathcal{R}^4$.
\end{cor}

\subsection{Some computations of $\sigma_{\mathcal{R}^4}$ and $\phi_{\mathcal{R}^4}$}
In this subsection we compute the value of our $\sigma_{\mathcal{R}^4}$ or $\phi_{\mathcal{R}^4}$
for some examples.

\begin{lem}
\label{lem:4-3-1}
Let $X=s_{3,3}(\mathbb{P}_1\times \mathbb{P}_1)$ be as defined in 4.2. Then $D_X$ is a hypersurface and
$\deg D_X=34$. For a lasso $\sigma_X$ around $D_X$, we have $\phi_X(\sigma_X)=-9/17$.
\end{lem}
\begin{proof}
Since $X\cong \mathbb{P}_1\times \mathbb{P}_1$, $X$ is simply connected and we have ${\rm Sign}X=0$, $\chi(X)=4$.
Also we have $\deg X=18$. By Proposition \ref{prop:3-5-2} and (\ref{eq:3-5-4}), the assertion follows.
\end{proof}

Let $\iota \colon \Delta \rightarrow \mathbb{P}(V_{3,3})$ be as in Proposition \ref{prop:2-3-3}.
Then we get a $\mathcal{R}^4$-degeneration $\iota^*\mathcal{W}\rightarrow \Delta$ (see Theorem \ref{thm:2-3-4}),
which we denote by $\mathcal{F}_I$ and call a \textit{singular fiber germ of type I}.
In this case we can choose $x_{\mathcal{F}_I}$ in (\ref{eq:4-2-2}) to be the image of a
lasso around $D_X$. By Lemma \ref{lem:3-5-1} the signature of the fiber neighborhood is $0$.
Thus we have

\begin{prop}
\label{prop:4-3-2}
$$\sigma_{\mathcal{R}^4}(\mathcal{F}_I)=\phi_{\mathcal{R}^4}(x_{\mathcal{F}_I})=-9/17.$$
\end{prop}

In the following let $\Delta=\{ z\in \mathbb{C};|z| \le \varepsilon \}$
for a sufficiently small real number $\varepsilon>0$.

\begin{exple}
\label{ex:4-1}
{\rm
Let
$$\Phi(z,a_0,a_1,b_0,b_1)=\varphi^0(a_0,a_1){b_0}^3
+({a_0}^3+z^6{a_1}^3)b_0{b_1}^2+z^9\varphi^3(a_0,a_1){b_1}^3$$
and $S_{\Phi} \subset \Delta \times \mathbb{P}_1\times \mathbb{P}_1$ the zero locus of $\Phi$.
Here $\varphi^0$, $\varphi^3$ are generic homogeneous polynomials of degree 3.
Let $f^{\prime}\colon S_{\Phi}\rightarrow \Delta$ be the first projection.

$S_{\Phi}$ has an isolated singularity at $(0,[0:1],[0:1])$.
Applying the resolution process given by Ashikaga \cite{Ashi}, we will obtain a resolution 
$\varpi\colon \widetilde{S}_{\Phi} \rightarrow S_{\Phi}$ of the singularity.
By successive blow down of $(-1)$-curves contained in the fiber at $0$,
we finally get the fiber germ what we want.

In the below we describe the resolution process.
We shall introduce the inhomogeneous coordinates $a=a_0/a_1$ and $b=b_0/b_1$.

First we recall a terminology from \cite{Ashi}.
Let $W$ be a complex manifold of dimension two and $L$ a holomorphic line bundle
on $W$. Let $\bar{L}=\mathbb{P}(\mathcal{O}_W\oplus \mathcal{O}_W(L))$
the associated $\mathbb{P}_1$ bundle. This is the $\mathbb{P}_1$ bundle on $W$
whose fiber at $w\in W$ is the projectivization of the dual space of $\mathbb{C}\oplus L_w$.
Here $L_w$ is the fiber of $L\rightarrow W$ at $w$.
Let $T=\mathcal{O}_{\bar{L}}(1)$. This is the line bundle on $\bar{L}$ whose
fiber at $\bar{\ell}$ (which is a line of $(\mathbb{C}\oplus L_w)^*$
for some $w\in W$) is the dual space of $\bar{\ell}$. Let $S$ be an irreducible
reduced divisor on $\bar{L}$ which is linearly equivalent to $3T$. In \cite{Ashi},
the triple $(S,W,L)$ is called a \textit{triple section surface}.

Let $W=\Delta \times \mathbb{P}_1$ with $(z,[a_0:a_1])$ the global coordinates, and $L$ a trivial line bundle
on $W$. We can introduce the homogeneous fiber coordinates $[b_0:b_1]$ for $\bar{L}$ by assigning
the linear functional on $\mathbb{C}\oplus L_w=\mathbb{C}\oplus \mathbb{C}$,
given by $(c_0,c_1)\mapsto b_0c_0+b_1c_1$, to $(b_0,b_1)$.
Our $S_{\Phi}$, which is the zero locus of $\Phi$, is naturally identified with an irreducible reduced
divisor on $\bar{L}$. Then $(S_{\Phi},W,L)$ is a triple section surface in the above sense.

Let $\tau_1\colon W_1\rightarrow W$ be the blow up at the origin $p_1=(0,[0:1])\in W$
and let $\hat{\tau}_1=(\bar{\tau_1},\tau_1)\colon (S_1,W_1,L_1)\rightarrow (S_{\Phi},W,L)$
be the triplet blow-up at $p_1$ with $\ell_1=1$, in the sense of \cite{Ashi}, p.181.
$W_1$ is covered by the two coordinate neighborhoods
$U_z=\{ (z,\tilde{a})\}$ and $U_a=\{ (\tilde{z},a)\}$, and $\tau_1$ is given by
$\tau_1(z,\tilde{a})=(z,\tilde{a}z)$ on $U_z$, and $\tau_1(\tilde{z},a)=(a\tilde{z},a)$ on $U_a$.
Note that $L_1=\tau_1^*L-E_1$ where $E_1$ is the exceptional curve of $\tau_1$.

Next, let $p_2=(0,0)\in U_z$ and $\tau_2\colon W_2\rightarrow W_1$ the blow up at $p_2$. Let
$\hat{\tau}_2=(\bar{\tau_2},\tau_2)\colon (S_2,W_2,L_2)\rightarrow (S_1,W_1,L_1)$
the triplet blow-up at $p_2$ with $\ell_2=2$.
$(S_2,W_2,L_2)$ is also a triple section surface.
Note that $L_2=\tau_2^*L_1-2E_2$ where $E_2$ is the exceptional curve of $\tau_2$.

There is a natural map $W_2\rightarrow W\rightarrow \Delta$, whose fiber at $0\in \Delta$ looks like Figure 1.
\begin{center}
\unitlength 0.1in
\begin{picture}( 20.0000, 12.6000)(  8.0000,-16.0000)
%
\special{pn 13}%
\special{pa 1000 400}%
\special{pa 1000 1400}%
\special{fp}%
%
\special{pn 13}%
\special{pa 800 600}%
\special{pa 2800 600}%
\special{fp}%
%
\special{pn 13}%
\special{pa 2400 400}%
\special{pa 2400 1600}%
\special{fp}%
\put(10.8000,-14.2000){\makebox(0,0)[lb]{$C$}}%
\put(14.7000,-5.1000){\makebox(0,0)[lb]{$N$}}%
\put(25.2000,-12.5000){\makebox(0,0)[lb]{$E_2$}}%
\end{picture}%

Figure 1
\end{center}
Here, $N$ is the proper transform of $E_1\subset W_1$, and $C$ is the proper
transform of $\{ z=0\} \subset W$. All the irreducible components are curves of genus 0.
We denote by $\hat{\pi}_2$ the natural projection $S_2\rightarrow W_2$ and
let $N^*={\hat{\pi}_2}^{-1}(N)$. $N^*\subset \bar{L}_2$ is a curve of genus 0.

Let $\bar{\sigma}\colon M\rightarrow \bar{L}_2$ be the blow up with center $N^*$
and let $\widetilde{S}_{\Phi}$ be the proper transform of $S_2$ by $\bar{\sigma}$
(see \cite{Ashi}, p.187). Then $\widetilde{S}_{\Phi}$ turns out to be non-singular.
Setting $\varpi$ to be the natural map from $\widetilde{S}_{\Phi}$ to $S_{\Phi}$,
we get a resolution $\varpi\colon \widetilde{S}_{\Phi}\rightarrow S_{\Phi}$.

The fiber $(f^{\prime}\circ \varpi)^{-1}(0)\subset \widetilde{S}_{\Phi}$ looks like Figure 2.
\begin{center}
\unitlength 0.1in
\begin{picture}( 24.0000, 18.0000)(  6.0000,-22.0000)
%
\special{pn 13}%
\special{pa 600 1000}%
\special{pa 1600 400}%
\special{fp}%
%
\special{pn 13}%
\special{pa 1000 600}%
\special{pa 1000 600}%
\special{fp}%
\special{pa 2600 600}%
\special{pa 2600 600}%
\special{fp}%
%
\special{pn 13}%
\special{pa 1000 600}%
\special{pa 1000 600}%
\special{fp}%
%
\special{pn 13}%
\special{pa 1000 600}%
\special{pa 3000 600}%
\special{fp}%
%
\special{pn 13}%
\special{pa 2600 400}%
\special{pa 2600 2000}%
\special{fp}%
%
\special{pn 13}%
\special{pa 3000 1400}%
\special{pa 1400 2200}%
\special{fp}%
%
\special{pn 13}%
\special{pa 2200 2000}%
\special{pa 600 2000}%
\special{fp}%
\put(17.2000,-8.7000){\makebox(0,0)[lb]{$N_1$}}%
\put(20.0000,-17.6000){\makebox(0,0)[lb]{$N_2$}}%
\put(26.8000,-11.6000){\makebox(0,0)[lb]{$\widetilde{E}_2$}}%
\put(9.1000,-11.4000){\makebox(0,0)[lb]{$C_1$}}%
\put(6.5000,-18.5000){\makebox(0,0)[lb]{$C_2$}}%
\end{picture}%

Figure 2
\end{center}
Here $C_1$, $N_1$, $N_2$ are curves of genus 0, $C_2$ is a curve
of genus 1, and $\widetilde{E}_2$ is a curve of genus 3. The self intersection numbers
are: ${C_1}^2={N_2}^2=-1$, ${N_1}^2={C_2}^2=-2$, and ${\widetilde{E}_2}^2=-3$.
The inverse images of $C$, $N$, and $E_2$ by the natural map $\widetilde{S}_{\Phi}\rightarrow W_2$
are $C_1\amalg C_2$, $N_1\amalg N_2$, and $\widetilde{E}_2$, respectively.

Note that we have a triple covering $\widetilde{S}_{\Phi}\rightarrow W_2$.
The restrictions of this map to $C_2$ or $\widetilde{E}_2$ gives
a double covering $C_2\rightarrow C$ with 4 simple branch points
or a triple covering $\widetilde{E}_2\rightarrow E_2$ with 10 simple branch points.
As a divisor, $(f^{\prime}\circ \varpi)^{-1}(0)=C_1+N_1+\widetilde{E}_2+2N_2+C_2$.

Finally let $\widetilde{S}_{\Phi}\rightarrow \bar{S}_{\Phi}$ be the contraction
obtained by repeating blow down of $(-1)$ curves in the fiber at $0\in \Delta$
until the resulting surface contains no more such curves (we need to blow down three times).
Let $\bar{C}_2$ (resp.\ $\bar{E}_2$) be the image of $C_2$ (resp.\ $\widetilde{E}_2$)
by this contraction. They are curves of genus 1 and 3 respectively, and
$\bar{C}_2\cdot \bar{E}_2=1$, ${\bar{C}_2}^2={\bar{E}_2}^2=-1$.

Let $f_{\Phi}\colon \bar{S}_{\Phi}\rightarrow \Delta$ be the map induced from $f^{\prime}\circ \varpi$.
Then $\mathcal{F}_{3,1}:=(\bar{S}_{\Phi},f_{\Phi},\Delta,0)$ is a $\mathcal{R}^4$-degeneration with
$f_{\Phi}^{-1}(0)$ being homeomorphic to the one point union of a surface of genus 1 and a surface of genus 3.
}
\end{exple}

\begin{prop}
\label{prop:4-3-3}
$$\sigma_{\mathcal{R}^4}(\mathcal{F}_{3,1})=11/17; \phi_{\mathcal{R}^4}(x_{\mathcal{F}_{3,1}})=28/17.$$
\end{prop}

\begin{proof}
The idea is to globalize the fiber germ $\mathcal{F}_{3,1}$ by a small perturbation.
As a result we will obtain a $\mathcal{R}^4$-fibration such that the set of singular fiber germs
consists of one $\mathcal{F}_{3,1}$ and $\mathcal{F}_I$'s.
Then we can compute the number of singular fiber germs, and get the value
of $\sigma_{\mathcal{R}^4}(\mathcal{F}_{3,1})$ by the global signature formula (\ref{eq:4-1-1}).

Let $\alpha \gg0$ be a sufficiently large natural number and let $\varphi=\varphi(a_0,a_1,b_0,b_1)$
be a generic $(3,3)$ homogeneous polynomial. Set $\Phi^{\prime}=\Phi+z^{\alpha}\varphi$.
Regarding $z$ as an affine coordinate of $\mathbb{P}_1$,
let $S_{\Phi^{\prime}}\subset \mathbb{P}_1\times \mathbb{P}_1\times \mathbb{P}_1$
be the zero locus of $\Phi^{\prime}$.

$S_{\Phi^{\prime}}$ has an isolated singularity at $(0,[0:1],[0:1])$,
and the same resolution process as that of $S_{\Phi}$ can be applied.
Let $\widetilde{S}_{\Phi^{\prime}}$ and $\bar{S}_{\Phi^{\prime}}$ be the
result of the process corresponding to $\widetilde{S}_{\Phi}$ and $\bar{S}_{\Phi}$,
respectively. Then $\widetilde{S}_{\Phi^{\prime}}$ is a non-singular compact complex
surface, and the induced projection $f_{\Phi^{\prime}}\colon \bar{S}_{\Phi^{\prime}}\rightarrow \mathbb{P}_1$
is a $\mathcal{R}^4$-fibration. The singular fiber germ at $0$ is $\mathcal{F}_{3,1}$,
and since $\psi$ is chosen to be generic, the other singular fiber germs is of type $I$.

We compute the various invariants.
First of all, the holomorphic Euler characteristic of $S_{\Phi^{\prime}}$ is computed as
$\chi(\mathcal{O}_{S_{\Phi^{\prime}}})=4\alpha-3$, and the self intersection number
of the dualizing sheaf (see Introduction) of $S_{\Phi^{\prime}}$ is computed as 
${\omega_{S_{\Phi^{\prime}}}}^2=14\alpha-24$.
By using Lemma 1.10 and the formula in p.187 of \cite{Ashi}, we have
$\chi(\mathcal{O}_{\widetilde{S}_{\Phi^{\prime}}})=4\alpha-10$ and
${\omega_{\widetilde{S}_{\Phi^{\prime}}}}^2=14\alpha-49$. Thus, we have
$\chi(\mathcal{O}_{\bar{S}_{\Phi^{\prime}}})=4\alpha-10$ and
${\omega_{\bar{S}_{\Phi^{\prime}}}}^2=14\alpha-46$.
By the Noether formula and the Hirzebruch signature formula, we have
$\chi(\bar{S}_{\Phi^{\prime}})=34\alpha-74$ and
${\rm Sign}(\bar{S}_{\Phi^{\prime}})=-18\alpha+34$.

As (\ref{eq:3-5-4}), the number of singular fiber germs of
$f_{\Phi^{\prime}}\colon \bar{S}_{\Phi^{\prime}}\rightarrow \mathbb{P}_1$ is computed as
$$\chi(\bar{S}_{\Phi^{\prime}})-2\cdot(2-2\cdot 4)=34\alpha-62.$$
In particular, the number of fiber germs of type I is $34\alpha-63$.
By the global signature formula, we have
$${\rm Sign}(\bar{S}_{\Phi^{\prime}})=-18\alpha+34=
(34\alpha-63)\cdot (-9/17)+\sigma_{\mathcal{R}^4}(\mathcal{F}_{3,1}).$$
Finally, since the signature of $\bar{S}_{\Phi}$ is $-1$,
$\phi_{\mathcal{R}^4}(\mathcal{F}_{3,1})=\sigma_{\mathcal{R}^4}(\mathcal{F}_{3,1})-1$.
This completes the proof.
\end{proof}

\begin{exple}
\label{ex:4-2}
{\rm
Let
$$\Phi(z,a_0,a_1,b_0,b_1)=\varphi^0(a_0,a_1){b_0}^3+\varphi^1(a_0,a_1){b_0}^2b_1+{a_0}^3b_0{b_1}^2
+z^6\varphi^3(a_0,a_1){b_1}^3$$
and $S_{\Phi}\subset \Delta \times \mathbb{P}_1 \times \mathbb{P}_1$ the zero locus of
$\Phi$. Here $\varphi^0$, $\varphi^1$, and $\varphi^3$ are generic homogeneous polynomials
of degree 3.

$S_{\Phi}$ has an isolated singularity at $(0,[0:1],[0:1])$.
If we write $\varphi^i(a)$ or $\Phi(z,a,b)$ instead of $\varphi^i(a,1)$ or $\Phi(z,a,1,b,1)$
respectively, then
$$\Phi(z,a,b)=(\varphi^0(a)b+\varphi^1(a))(b^2+A^3b+Z^6),$$
where $A=a(\varphi^0(a)b+\varphi^1(a))^{-1/3}$ and
$Z=z\{\varphi^3(a)/(\varphi^0(a)b+\varphi^1(a))\}^{1/6}$.
Note that
$$b^2+a^3b+z^6=\left(b+\frac{a^3}{2}\right)^2-\frac{a^6}{4}+Z^6.$$
Set $h(b,z,a)=b^2+z^6-a^6$ and let $S_h \subset \mathbb{C}^3$ be the zero locus of $h$.
Define $f^{\prime}\colon S_h\rightarrow \Delta$ by $f^{\prime}(b,z,a)=z$.
Then the singularity of $S_{\Phi}$ is analytically equivalent to the
hypersurface singularity $(S_h,0)$.
In this case Horikawa's canonical resolution for double coverings (see \cite{Hori}, section 2)
can be used for a resolution of $(S_h,0)$.

The process is as follows. Let $W$ be an open neighborhood of $0$ in $\mathbb{C}^2$
with global coordinates $(z,a)$ and
let $\tau\colon \widetilde{W}\rightarrow W$ be the blow up at the origin.
$\widetilde{W}$ is covered by the two coordinate neighborhoods $U_z$ and $U_a$ as in Example \ref{ex:4-1},
and $\tau$ is given by the same formula as $\tau_1$.
The picture of $\widetilde{W}$ is like Figure 3.
\begin{center}
\unitlength 0.1in
\begin{picture}( 24.4000, 10.2000)(  5.6000,-12.0000)
%
\special{pn 13}%
\special{pa 1000 1200}%
\special{pa 1000 400}%
\special{fp}%
\special{sh 1}%
\special{pa 1000 400}%
\special{pa 980 468}%
\special{pa 1000 454}%
\special{pa 1020 468}%
\special{pa 1000 400}%
\special{fp}%
%
\special{pn 13}%
\special{pa 2600 1200}%
\special{pa 2600 400}%
\special{fp}%
\special{sh 1}%
\special{pa 2600 400}%
\special{pa 2580 468}%
\special{pa 2600 454}%
\special{pa 2620 468}%
\special{pa 2600 400}%
\special{fp}%
%
\special{pn 13}%
\special{pa 1800 800}%
\special{pa 3000 800}%
\special{fp}%
\special{sh 1}%
\special{pa 3000 800}%
\special{pa 2934 780}%
\special{pa 2948 800}%
\special{pa 2934 820}%
\special{pa 3000 800}%
\special{fp}%
%
\special{pn 13}%
\special{pa 1800 800}%
\special{pa 600 800}%
\special{fp}%
\special{sh 1}%
\special{pa 600 800}%
\special{pa 668 820}%
\special{pa 654 800}%
\special{pa 668 780}%
\special{pa 600 800}%
\special{fp}%
\put(10.9000,-4.2000){\makebox(0,0)[lb]{$z$}}%
\put(5.6000,-6.9000){\makebox(0,0)[lb]{$\tilde{a}$}}%
\put(26.4000,-3.5000){\makebox(0,0)[lb]{$a$}}%
\put(28.6000,-10.4000){\makebox(0,0)[lb]{$\tilde{z}$}}%
\put(23.8000,-11.9000){\makebox(0,0)[lb]{$C$}}%
\put(18.4000,-7.2000){\makebox(0,0)[lb]{$E$}}%
\end{picture}%

Figure 3
\end{center}
Here $E$ is the exceptional curve of $\tau$ and $C$ is the proper transform of $\{z=0\}\subset W$.

Let $\widetilde{S}_h$ be the double covering on $\widetilde{W}$ defined by the following:
on $U_z$, let $\widetilde{S}_h^z$ be the hypersurface in $\mathbb{C}^3$ given by
$x^2+1-{\tilde{a}}^6=0$ and on $U_a$, let $\widetilde{S}_h^a$ be given by $x^2+{\tilde{z}}^6-1=0$.
Identifying $\widetilde{S}_h^z$ with $\widetilde{S}_h^a$ over $U_z\cap U_a$ by the map
$(x,z,\tilde{a})\mapsto (x/{\tilde{a}}^3,\tilde{z},a)$, we get the resulting
non-singular surface $\widetilde{S}_h$.
Define $\varpi\colon \widetilde{S}_h\rightarrow S_h$ by
$\varpi(x,z,\tilde{a})=(xz^3,z,\tilde{a})$ on $\widetilde{S}_h^z$, and
$\varpi(x,\tilde{z},a)=(xa^3,a\tilde{z},a)$ on $\widetilde{S}_h^a$. Then this is a resolution of singularity.
$(f^{\prime}\circ \varpi)^{-1}(0)$ looks like:

\begin{center}
\unitlength 0.1in
\begin{picture}( 18.0000,  6.4000)(  4.0000,-10.0000)
%
\special{pn 13}%
\special{pa 400 600}%
\special{pa 2200 600}%
\special{fp}%
%
\special{pn 13}%
\special{pa 1600 400}%
\special{pa 1600 1000}%
\special{fp}%
%
\special{pn 13}%
\special{pa 2000 400}%
\special{pa 2000 1000}%
\special{fp}%
\put(6.5000,-5.3000){\makebox(0,0)[lb]{$\widetilde{E}$}}%
\put(13.3000,-9.9000){\makebox(0,0)[lb]{$C_1$}}%
\put(20.5000,-9.6000){\makebox(0,0)[lb]{$C_2$}}%
\end{picture}%

Figure 4
\end{center}
Here $\widetilde{E}$ (resp.\ $C_1\amalg C_2$) is the proper transform of $E$ (resp.\ $C$).
$\widetilde{E}$ is a curve of genus 2.
Note that the restriction of $\widetilde{S}_h \rightarrow \widetilde{W}$ to $\widetilde{E}$
gives a double covering $\widetilde{E}\rightarrow E\cong \mathbb{P}_1$ with 6 simple branch points.

Now applying the above process, we can resolve the singularity
of $S_{\Phi}$. Let $\varpi\colon \widetilde{S}_{\Phi}\rightarrow S_{\Phi}$ be such a resolution.
Then $(f^{\prime}\circ \varpi)^{-1}(0)$ consists of the three components, $\widetilde{E}$,
$C_1$, and $C_2$, like Figure 4 (using the same letter). Note that
$$\Phi(0,a,b)=b(\varphi^0(a)b^2+\varphi^1(a)b+a^3).$$
We assume that $C_1$ (resp.\ $C_2$) is the component corresponding to $b=0$
(resp.\ $\varphi^0(a)b^2+\varphi^1(a)b+a^3=0$).
$C_1$ is a curve of genus 0 and $C_2$ is a curve of genus 2.
The self intersection numbers are: ${\widetilde{E}}^2=-2$, ${C_1}^2=-1$, ${C_2}^2=-1$.

Let $\widetilde{S}_{\Phi}\rightarrow \bar{S}_{\Phi}$ be the blow down of $\bar{C}_1$
and $f_{\Phi}\colon \bar{S}_{\Phi}\rightarrow \Delta$ the map induced from $f^{\prime}\circ \varpi$.
Then $\mathcal{F}_{2,2}:=(\bar{S}_{\Phi},f_{\Phi},\Delta,0)$ is a $\mathcal{R}^4$-degeneration with
$f_{\Phi}^{-1}(0)$ being homeomorphic to the one points union of two curves of genus 2.
}
\end{exple}

\begin{prop}
\label{prop:4-3-4}
$$\sigma_{\mathcal{R}^4}(\mathcal{F}_{2,2})=19/17; \phi_{\mathcal{R}^4}(x_{\mathcal{F}_{2,2}})=36/17.$$
\end{prop}

\begin{proof}
The idea of proof is the same as Proposition \ref{prop:4-3-3}.
Let $\alpha$, $\varphi$, $\Phi^{\prime}$, and $S_{\Phi^{\prime}}$
be the same as in the proof of Proposition \ref{prop:4-3-3}. We construct $\widetilde{S}_{\Phi^{\prime}}$ and
$f_{\Phi^{\prime}}\colon \bar{S}_{\Phi^{\prime}}\rightarrow \mathbb{P}_1$
by a similar manner to Proposition \ref{prop:4-3-3} except for using the resolution of $(S_h,0)$ described as above.

We have $\chi(\mathcal{O}_{S_{\Phi^{\prime}}})=4\alpha-3$, and
${\omega_{S_{\Phi^{\prime}}}}^2=14\alpha-24$.
Using Lemma 6 of \cite{Hori}, we have
$\chi(\mathcal{O}_{\widetilde{S}_{\Phi^{\prime}}})=4\alpha-6$, and
${\omega_{\widetilde{S}_{\Phi^{\prime}}}}^2=14\alpha-32$. Thus we have
$\chi(\mathcal{O}_{\bar{S}_{\Phi^{\prime}}})=4\alpha-6$ and
${\omega_{\bar{S}_{\Phi^{\prime}}}}^2=14\alpha-31$, therefore
$\chi(\bar{S}_{\Phi^{\prime}})=34\alpha-41$ and ${\rm Sign}(\bar{S}_{\Phi^{\prime}})=-18\alpha+17$.
Now the number of singular fiber germs of $f_{\Phi^{\prime}}$ is
$$34\alpha-41-2\cdot (2-2\cdot 4)=34\alpha-29,$$
hence the number of singular fiber germs of type I is $34\alpha-30$.
By the global signature formula we have
$$-18\alpha+17=(34\alpha-30)\cdot (-9/17)+\sigma_{\mathcal{R}^4}(\mathcal{F}_{2,2}).$$
Finally the signature of $\bar{S}_{\Phi}$ is $-1$. This completes the proof.
\end{proof}

\begin{exple}
\label{ex:4-3}
{\rm
Let
$$\Phi(z,a_0,a_1,b_0,b_1)=(a_0b_1-a_1b_0)^3+z^3\varphi^0(a_0,a_1,b_0,b_1).$$
We write $Y=\Delta \times \mathbb{P}_1 \times \mathbb{P}_1$ and let
$S_{\Phi}\subset Y$ be the zero locus of $\Phi$.
Here $\varphi^0\in V_{3,3}$ is a generic $(3,3)$ homogeneous polynomial.

Let $\Gamma=\{z=0\} \subset Y$ and we denote by $D$ the fiber at $0$ of the first projection
$S_{\Phi}\rightarrow \Delta$. $D$ is the diagonal locus in $\Gamma\cong \mathbb{P}_1\times \mathbb{P}_1$.
Let $\tau_1\colon Y_1\rightarrow Y$ be the blow up along $D$.
Let $E_0\subset Y_1$ be the proper transform of $\Gamma$, $E_1$ the exceptional set of $\tau_1$,
and $\widetilde{S}_{\Phi}\subset Y_1$ the proper transform of $S_{\Phi}$.
Note that $E_1$ is isomorphic to the Hirzebruch surface $\mathbb{F}_2$ of degree 2.

Then we see that $\widetilde{S}_{\Phi}$ is non-singular and $\widetilde{S}_{\Phi}\cap E_0=\emptyset$.
We write $f_{\Phi}\colon \widetilde{S}_{\Phi}\rightarrow \Delta$ the natural projection.
Then we have a $\mathcal{R}^4$-degeneration $\mathcal{F}_R^{\prime}:=(\widetilde{S}_{\Phi},f_{\Phi},\Delta,0)$.
We see that $f_{\Phi}^{-1}(0)=E_1\cap S_1$ is a smooth curve of genus 4.
This curve is non-hyperelliptic, but not a curve of rank 4. This can be seen as follows.
First we can contract $E_0$. Let $\bar{\tau}\colon Y_1\rightarrow \bar{Y}$ be the contraction of $E_0$.
The projection $Y_1\rightarrow \Delta$ induces the projection $\tilde{f}\colon \bar{Y}\rightarrow \Delta$,
whose central fiber ${\tilde{f}}^{-1}(0)=\bar{\tau}(E_1)$ is $\mathbb{F}_2$ and
the other fibers are isomorphic to $\mathbb{P}_1\times \mathbb{P}_1$. Thus, we can think
$f_{\Phi}^{-1}(0)$ is contained in $\mathbb{F}_2$.
Contracting the negative section of $\mathbb{F}_2$, we get a quadric $Q_3$ of rank 3
in $\mathbb{P}_3$. If we map $f_{\Phi}^{-1}(0)$ into $\mathbb{P}_3$ by this contraction, then
$f_{\Phi}^{-1}(0)$ can be realized as a $(2,3)$-complete intersection: the intersection
of $Q_3$ and some cubic surface. Thus $f_{\Phi}^{-1}(0)$ is non-hyperelliptic but not of rank 4.
Topologically $f_{\Phi}\colon \widetilde{S}_{\Phi}\rightarrow \Delta$ is a trivial $\Sigma_4$-bundle.
}
\end{exple}

\begin{prop}
\label{prop:4-3-5}
Let $\mathcal{F}_R^{\prime}=(\widetilde{S}_{\Phi},f,\Delta,0)$ be the fiber germ as above. Then we have
$\sigma_{\mathcal{R}^4}(\mathcal{F}_R^{\prime})=\phi_{\mathcal{R}^4}(x_{\mathcal{F}_R^{\prime}})=4/17$.
\end{prop}

\begin{proof}
The proof proceeds as the same before. Let $\alpha$, $\varphi$, $\Phi^{\prime}$,
$S_{\Phi^{\prime}}$, and $\widetilde{S}_{\Phi^{\prime}}$
be the same as in the proof of Proposition \ref{prop:4-3-3}.
We use the same notation $Y$ and $Y_1$ for $\mathbb{P}_1\times \mathbb{P}_1\times \mathbb{P}_1$ and
the blow up along $D$, respectively.

First we have $\chi(\mathcal{O}_{S_{\Phi^{\prime}}})=4\alpha-3$ and
${\omega_{S_{\Phi^{\prime}}}}^2=14\alpha-24$. We need to compute
$\chi(\mathcal{O}_{\widetilde{S}_{\Phi^{\prime}}})$ and
${\omega_{\widetilde{S}_{\Phi^{\prime}}}}^2$.
By $K_{Y_1}\sim \tau_1^*K_Y+E_1$ (linear equivalence) and $S_1\sim \tau_1^*S-3E_1$, we have
$K_{Y_1}+S_1\sim \tau_1^*(K_Y+S)-2E_1$. By restricting to $S_1$, we have
$\tau_1^*\omega_S=\omega_{S_1}+2E_1|_{S_1}$. We get
\begin{eqnarray*}
{\omega_S}^2=(\tau_1^*\omega_S)^2 &=& (\omega_{S_1}+2E_1|_{S_1})^2 \\
&=& {\omega_{S_1}}^2+4\omega_{S_1}\cdot E_1|_{S_1}+4(E_1|_{S_1})^2={\omega_{S_1}}^2+24,
\end{eqnarray*}
since $\omega_{S_1}\cdot E_1|_{S_1}=2g(E_1\cap S_1)-2-(E_1|_{S_1})^2=6-(E_1|_{S_1})^2$ by the adjunction formula.
Therefore, ${\omega_{\widetilde{S}_{\Phi^{\prime}}}}^2=14\alpha-48$. We next compute
\begin{eqnarray*}
\chi(\mathcal{O}_{S_1}) &=& \chi(\mathcal{O}_{Y_1})-\chi(\mathcal{O}_{Y_1}(-S_1)) \\
&=& \chi(\mathcal{O}_{Y_1})-\chi(\mathcal{O}_{Y_1}(-S_1-E_1))-\chi(\mathcal{O}_{E_1}(-S_1|_{E_1})) \\
&=& \cdots \\
&=& \chi(\mathcal{O}_{Y_1})-\chi(\mathcal{O}_{Y_1}(-S_1-3E_1))
-\sum_{i=0}^2\chi(\mathcal{O}_{E_1}(-S_1|_{E_1}-iE_1|_{E_1})).
\end{eqnarray*}
We have $\chi(\mathcal{O}_{Y_1})-\chi(\mathcal{O}_{Y_1}(-S_1-3E_1))
=\chi(\mathcal{O}_{Y_1})-\chi(\mathcal{O}_{Y_1}(-\tau_1^*S))=\chi(\mathcal{O}_Y)-\chi(\mathcal{O}_Y(-S))
=\chi(\mathcal{O}_{S})$. To compute the remaining term, we use
$$\chi(\mathcal{O}_{E_1}(-S_1|_{E_1}-iE_1|_{E_1}))=\chi(\mathcal{O}_{E_1}(-iE_1|_{E_1}))-
\chi(\mathcal{O}_{E_1\cap S_1}(-iE_1|_{E_1\cap S_1})).$$

Note that the divisor $(E_0+E_1)|_{E_1}$ is trivial on $E_1$,
since $E_0+E_1$ is a fiber of $Y_1\rightarrow \Delta$;
$C_{\infty}=E_0\cap E_1$ is the negative section of $E_1\cong \mathbb{F}_2$, so ${C_{\infty}}^2=-2$;
$E_0\cap E_1\cap S_1=\emptyset$. From these we have
$\mathcal{O}_{E_1}(-E_1|_{E_1})\cong \mathcal{O}_{E_1}(C_{\infty})$ and
$\mathcal{O}_{E_1\cap S_1}(-E_1)\cong \mathcal{O}_{E_1\cap S_1}$. Using the Riemann-Roch formula
and $\chi(\mathcal{O}_{E_1})=1$, we get $\chi(\mathcal{O}_{E_1}(-S_1|E_1))=4$;
$\chi(\mathcal{O}_{E_1}(-S_1|_{E_1}-E_1|_{E_1}))=3$; $\chi(\mathcal{O}_{E_1}(-S_1|_{E_1}-2E_1|_{E_1}))=0$.

In summary, we have $\chi(\mathcal{O}_{S_1})=\chi(\mathcal{O}_{S})-7$. Therefore,
$\chi(\mathcal{O}_{\widetilde{S}_{\Phi^{\prime}}})=4\alpha-10$.

Now $\chi(\widetilde{S}_{\Phi^{\prime}})=34\alpha-72$, and
${\rm Sign}(\widetilde{S}_{\Phi^{\prime}})=-18\alpha+32$. The number of topologically singular fibers is
$$34\alpha-72-2\cdot (2-2\cdot 4)=34\alpha-60.$$
By the global signature formula we have
$$-18\alpha+32=(34\alpha-60)\cdot (-9/17)+\sigma_{\mathcal{R}^4}(\mathcal{R}_3).$$
Finally, the signature of $\widetilde{S}_{\Phi}$ is zero. This completes the proof.
\end{proof}

\begin{exple}
\label{ex:4-4}
{\rm
In this last example we do not use the global signature
formula (\ref{eq:4-1-1}) directly.
Let $q_1={x_0}^2+{x_1}^2+{x_3}^2$ and $q_2={x_1}^2+{x_2}^2-{x_3}^2$, and
$h(x_0,x_1,x_2,x_3)$ a generic cubic polynomial. Let
$$S_h=\{ (x,z)\in \mathbb{P}_3\times \Delta; h(x)=(q_1+zq_2)(x)=0\},$$
and let $f\colon S_h\rightarrow \Delta$ be the second projection.
Note that $q_1+zq_2$ defines a smooth quardric except for $z=0$ and $f^{-1}(0)$ is contained in
the singular quadric $\{ q_1=0\}$. Thus $\mathcal{F}_R:=(S_h,f,\Delta,0)$ is a $\mathcal{R}^4$-degeneration.
Topologically $f\colon S_h\rightarrow \Delta$ is a trivial $\Sigma_4$-bundle.
}
\end{exple}

\begin{prop}
\label{prop:4-3-6}
$$\sigma_{\mathcal{R}^4}(\mathcal{F}_R)=\phi_{\mathcal{R}^4}(x_{\mathcal{F}_R})=2/17.$$
\end{prop}

\begin{proof}
We start from describing the associated principal $\mathcal{G}$-bundle $P(\xi_{\mathcal{F}_R})$.

Let $\omega_z$ be the basis of $\Omega^1(f^{-1}(z))\subset \mathbb{P}_3$ corresponding
to the homogeneous coordinates of $\mathbb{P}_3$ (see the paragraph before Theorem \ref{thm:4-2-1}).
Any frame of $\Omega^1(f^{-1}(z))$ (modulo $\mathbb{C}^*$) is written as the form $A\cdot \omega_z$,
$A\in PGL(4)$. Then
$$P(\xi_{\mathcal{F}_R})\cong \left\{ (z,A)\in (\Delta \setminus \{0\})\times PGL(4);
(A^{-1})^t B(z) A^{-1}=H \right\},$$
where
$$B(z)=\left( \begin{array}{cccc}
1 & 0 & 0 & 0 \\
0 & 1+z & 0 & 0 \\
0 & 0 & z & 0 \\
0 & 0 & 0 & 1-z \\
\end{array} \right).$$
Let $\kappa\colon \Delta \setminus \{0\}\rightarrow \Delta \setminus \{0\}$ be the map defined by $w\mapsto z^2$
and consider the pull back $\kappa^*P(\xi_{\mathcal{F}_R})$. Then the principal $\mathcal{G}$-bundle
$\kappa^*P(\xi_{\mathcal{F}_R})\rightarrow \Delta \setminus \{0\}$ has a section given by
$w\mapsto (w,A(w))$ where
$$A(w)=\frac{1}{\sqrt{2}} \left( \begin{array}{cccc}
1 & \sqrt{-1}r_1(w) & 0 & 0 \\
0 & 0 & w & \sqrt{-1}r_2(w) \\
0 & 0 & -w & \sqrt{-1}r_2(w) \\
1 & -\sqrt{-1}r_1(w) & 0 & 0 \\
\end{array} \right).$$
Here, $r_1(w)=\sqrt{1+w^2}$ and $r_2(w)=\sqrt{1-w^2}$. Thus, as a candidate for
$g_{\kappa^*\mathcal{F}_R}$ (see the proof of Proposition \ref{prop:4-1-7}), we can take
the map given by
$w\mapsto [e,A(w)\cdot f(a_0b_0,a_0b_1,a_1b_0,a_1b_1)]$. Since the diagram
$$\xymatrix{
 & \pi_1(U^X_{\mathcal{G}}) \\
 \pi_1(\Delta \setminus \{0\}) \ar[r]^{\times 2} \ar[ur]^{g_{\kappa^*\mathcal{F}_R}}
& \pi_1(\Delta \setminus \{0\}) \ar[u]_{g_{\mathcal{F}_R}} \\}$$
is commutative up to conjugacy,
${x_{\mathcal{F}_R}}^2$ is conjugate to
$x_{\kappa^*\mathcal{F}_R}=g_{\kappa^*\mathcal{F}_R}(\partial \Delta)$.

But now $A(w)\cdot f(a_0b_0,a_0b_1,a_1b_0,a_1b_1)$ is equal to
$$f\left(\frac{a_0b_0+a_1b_1}{\sqrt{2}},\frac{a_0b_0-a_1b_1}{\sqrt{-1}r_1(w)},
\frac{a_0b_1-a_1b_0}{w},\frac{a_0b_1+a_1b_0}{\sqrt{-1}r_2(w)}\right).$$
Modulo $\mathbb{C}^*$ this can be written as
$$(a_0b_1-a_1b_0)^3+w^3(\varphi^0+w\varphi^1+{\rm higher\ term\ with\ respect\ to\ }w)$$
where $\varphi^i$ is some $(3,3)$ homogeneous polynomial. This shows that $x_{\kappa^*\mathcal{F}_R}$
is homotopic to $x_{\mathcal{F}_R^{\prime}}$ in Proposition \ref{prop:4-3-5}.
Since $f\colon S_h\rightarrow \Delta$ is topologically trivial,
$$4/17=\phi_{\mathcal{R}^4}(x_{\mathcal{F}_R^{\prime}})=
\phi_{\mathcal{R}^4}(x_{\kappa^*\mathcal{F}_R})
=\phi_{\mathcal{R}^4}({x_{\mathcal{F}_R}}^2)=2\phi_{\mathcal{R}^4}(x_{\mathcal{F}_R}).$$
This completes the proof.
\end{proof}

Compare the above proof with \cite{AY}, Example 7.5, where the same fiber germ is considered.

\subsection{Fibrations of non-trigonal curves of genus 5}
Let $C$ be a non-hyperelliptic Riemann surface of genus 5. By the Enriques-Petri theorem (\cite{GH}, p.535),
the canonical image of $C$ is a $(2,2,2)$ complete intersection iff $C$ is non-trigonal, i.e.,
there is no holomorphic map $C\rightarrow \mathbb{P}_1$ of degree 3. Let $\mathcal{NT}^5\subset \mathbb{M}_5$
be the set of non-hyperelliptic and non-trigonal Riemann surface of genus 5. $\mathcal{NT}^5$ is Zariski
open in $\mathbb{M}_5$.

We denote by $[\alpha_0:\alpha_1:\alpha_2:\alpha_3:\alpha_4]$
the homogeneous coordinates of $\mathbb{P}_4$ and let $S_5$ be the space of $5\times 5$ symmetric matrices.
The Veronese map $v_2\colon \mathbb{P}_4\rightarrow \mathbb{P}(S_5)$ is given by
$$v_2([\alpha_0:\alpha_1:\alpha_2:\alpha_3:\alpha_4]):=[(\alpha_0,\alpha_1,\alpha_2,\alpha_3,\alpha_4)^t \cdot
(\alpha_0,\alpha_1,\alpha_2,\alpha_3,\alpha_4)].$$
This map is equivariant with respect to the action of $\mathcal{G}=PGL(5)$, where the action of $\mathcal{G}$
on $\mathbb{P}_4$ is induced by
$$A\cdot (\alpha_0,\alpha_1,\alpha_2,\alpha_3,\alpha_4)=(\alpha_0,\alpha_1,\alpha_2,\alpha_3,\alpha_4)A^t$$
for $A\in GL_5(\mathbb{C})$ and $\alpha_i \in \mathbb{C}$, and the action of $\mathcal{G}$ on $\mathbb{P}(S_5)$ is
induced by $A\cdot B=ABA^t$ for $A\in GL_5(\mathbb{C})$ and $B\in S_5$. Set
$$X:=v_2(\mathbb{P}_4).$$
The action of $\mathcal{G}$ on $G_{11}(\mathbb{P}(S_5))$ induced from the $\mathcal{G}$-action on $\mathbb{P}(S_5)$
preserves $D_X$ and $U^X=G_{11}(\mathbb{P}(S_5)) \setminus D_X$, and the projection
$p_X\colon \mathcal{C}^X\rightarrow U^X$ is $\mathcal{G}$-equivariant. Note that for $W\in U^X$, the fiber
$p_X^{-1}(W)$ is isomorphic to a smooth complete intersection in $\mathbb{P}_4$ of type $(2,2,2)$.
Thus $\xi^X=(\mathcal{C}^X,p_X,U^X)$ is a $\mathcal{NT}^5$-family.

We claim that $\xi^X$ and the $\mathcal{G}$-action on it satisfies the conditions in Proposition \ref{prop:4-1-3}.
The proof is similar to the case of $\mathcal{R}^4$, so we only describe the way of constructing
$\mathcal{G}$-bundles and $\mathcal{G}$-equivariant maps.

Let $C\in \mathcal{NT}^5$. By Max Noether's theorem, the natural map
$t_2\colon {\rm Sym}^2\Omega^1(C)\rightarrow H^0(C;K_C^{\otimes 2})$ is surjective hence the kernel is
3-dimensional. By taking the dual, we get the codimension 3 subspace
${\rm Ker}(t_2)^* \subset {\rm Sym}^2\Omega^1(C)^*$.
If we take a basis $\omega=(\omega_0,\omega_1,\omega_2,\omega_3,\omega_4)$ of $\Omega^1(C)$,
${\rm Sym}^2\Omega^1(C)^*$ is identified with $S_5$ by assigning $B\in S_5$ with the quadratic function
$$\Omega^1(C)\rightarrow \mathbb{C}, \ x_0\omega_0+\cdots +x_4\omega_4 \mapsto
(x_0,\ldots,x_4)B(x_0,\ldots,x_4)^t,$$
and we have the plane $B(\omega)\subset \mathbb{P}(S_5)$ of codimension 3 corresponding to ${\rm Ker}(t_2)^*$.
Note that the image $v_2\circ \iota_{\omega}(C)$, where $\iota_{\omega}\colon C\rightarrow \mathbb{P}_4$
is the canonical map given by $c\mapsto [\omega_0(c):\ldots:\omega_4(c)]$,
is the smooth intersection of $X$ and $B(\omega)$.

Let $P(C)$ be the set of all frames of $\Omega^1(C)$ modulo $\mathbb{C}^*$. By assigning $\omega \in P(C)$
with $B(\omega)$, we have the map $E_C\colon P(C)\rightarrow U^X$. Moreover, given the left action of
$\mathcal{G}$ on $P(C)$ by the same formula as the action of $PO_4^H(\mathbb{C})$ on $P(C)$ in subsection 4.2,
we see that $E_C$ is $\mathcal{G}$-equivariant.

Let $\xi=(\mathcal{C},p,B)$ be a $\mathcal{NT}^5$-family.
Applying the above construction to all the fibers,
we get a principal $\mathcal{G}$-bundle $P(\xi)$ and a $\mathcal{G}$-equivariant map
$E_{\xi}\colon P(\xi)\rightarrow U^X$ what we want.

\begin{thm}
\label{thm:4-4-1}
Let $\mathcal{NT}^5$ be the set of non-hyperelliptic and non-trigonal Riemann surfaces of genus 5 and $X$,
$\mathcal{G}$ as above. Then $\xi^X_{\mathcal{G}}=(\mathcal{C}^X_{\mathcal{G}},p_u,U^X_{\mathcal{G}})$
is a universal $\mathcal{NT}^5$-family.
We denote by $\rho_u\colon \pi_1(U^X_{\mathcal{G}})\rightarrow \Gamma_5$
the topological monodromy of $p_u\colon \mathcal{C}^X_{\mathcal{G}}\rightarrow U^X_{\mathcal{G}}$.
Then there exists a unique $\mathbb{Q}$-valued 1-cochain
$\phi_{\mathcal{NT}^5}\colon \pi_1(U^X_{\mathcal{G}})\rightarrow \mathbb{Q}$ whose coboundary
equals to $\rho_u^*\tau_5$.
\end{thm}
\begin{proof}
The proof is the same as the proof of Theorem \ref{thm:4-2-1} except that (\ref{eq:4-2-1}) is replaced
with the exact sequence
$$\pi_1(\mathcal{G})\cong \mathbb{Z}/5\mathbb{Z} \rightarrow \pi_1(U^X)\rightarrow \pi_1(U^X_{\mathcal{G}})
\rightarrow *.$$
\end{proof}

\begin{cor}
\label{cor:4-4-2}
Let $\mathcal{NT}^5$ be the set of non-hyperelliptic and non-trigonal Riemann surfaces of genus 5.
Then the formula
$$\sigma_{\mathcal{NT}^5}(\mathcal{F})=\phi_{\mathcal{NT}^5}(x_{\mathcal{F}})+{\rm Sign}(S)$$
for $\mathcal{F}=(S,\pi,\Delta,b)\in \mathcal{NT}^5_{loc}$ (see (\ref{eq:4-1-2})) gives
a local signature with respect to $\mathcal{NT}^5$.
\end{cor}

\begin{lem}
\label{lem:4-4-3}
Let $X=v_2(\mathbb{P}_4)$ be as above. Then $D_X$ is a hypersurface and
$\deg D_X=40$. For a lasso $\sigma_X$ around $D_X$, we have $\phi_X(\sigma_X)=-1/2$.
\end{lem}
\begin{proof}
This follows from Proposition \ref{prop:3-5-4}.
\end{proof}

Let $\iota\colon \Delta \rightarrow G_{11}(\mathbb{P}(S_5))$ be as in Proposition \ref{prop:2-3-3}.
Then we get a $\mathcal{NT}^5$-degeneration $\iota^*\mathcal{W}\rightarrow \Delta$, which
we denote by $\mathcal{F}_I$ and call a singular fiber germ of type I.

\begin{prop}
\label{prop:4-4-4}
Let $\mathcal{F}_I\in \mathcal{NT}_{loc}^5$ be as above. Then
$$\sigma_{\mathcal{NT}^5}(\mathcal{F}_I)=\phi_{\mathcal{NT}^5}(x_{\mathcal{F}_I})=-1/2.$$
\end{prop}
\begin{proof}
The proof is similar to the proof of Proposition \ref{prop:4-3-2}.
\end{proof}

This fiber germ is expected to play an important role when computing examples
as $\mathcal{F}_I \in \mathcal{R}^4_{loc}$ behaves like an "atomic" germ in subsection 4.3,
but at the present moment we don't have any example of element of $\mathcal{NT}^5$
other than $\mathcal{F}_I$ whose local signature has been computed.

\vspace{0.5cm}
\noindent \textbf{Final remarks.}
Although the construction of our local signature is purely topological,
we have used some algebraic geometry to compute examples.
It is an interesting problem to find and compute examples of fiber germs beyond the reach of
algebraic geometry, or to give a formula for the Meyer functions $\phi_X$ as Meyer and Atiyah did.
To do this we need a greater understanding
of $\rho_X$ or the topological monodromy $\rho_u$ of a universal $\mathcal{A}$-family.

In the case of $\mathcal{A}=\mathcal{R}^4$ or $\mathcal{NT}^5$, $\mathcal{A}$ is Zariski open in the moduli
space. Using this, we can prove that $\rho_u$ is surjective.
The proof is similar to the proof of \cite{Ku}, Proposition 6.3. Here is an outline.
Let $\mathcal{T}_g$ be the Teichm\"uller space of genus $g$,
and let $\tilde{\mathcal{A}}\subset \mathcal{T}_g$
be the inverse image of $\mathcal{A}$ by the quotient map $\mathcal{T}_g\rightarrow \mathbb{M}_g$.
Then $\tilde{\mathcal{A}}$ is Zariski open hence path connected, and is preserved by the action
of $\Gamma_g$ on $\mathcal{T}_g$.
By a natural way we get an $\mathcal{A}$-family on the Borel construction
$\tilde{\mathcal{A}}_{\Gamma_g}$, which is easily seen to be universal.
The homotopy exact sequence
$$\pi_1(\tilde{\mathcal{A}}_{\Gamma_g})\rightarrow \pi_1(B\Gamma_g)=\Gamma_g \rightarrow
\pi_0(\tilde{\mathcal{A}})=*$$
of the $\tilde{\mathcal{A}}$-bundle $\tilde{\mathcal{A}}_{\Gamma_g}\rightarrow B\Gamma_g$
shows the desired surjectivity.

\vspace{0.5cm}
\noindent \textbf{Acknowledgements.} I would like to thank Tadashi Ashikaga,
who kindly communicated the constructions of the first three examples in subsection 4.3 to me.
I also would like to thank my advisor Nariya Kawazumi for reading a draft, giving various comments,
and a warm encouragement. This research is supported by JSPS Research Fellowships for Young
Scientists (19$\cdot$5472).

\noindent \textsc{Yusuke Kuno\\
Department of Mathematics,
Graduate School of Science,\\
Hiroshima University,\\
Kagamiyama, Higashi-Hiroshima, 739-8526,
JAPAN}

\noindent \texttt{E-mail address:
kunotti@hiroshima-u.ac.jp}

\end{document}